\documentclass[10pt]{amsart}
\usepackage{amsmath}
\usepackage{amssymb}
\usepackage[final]{graphicx}
\usepackage{float}
\usepackage{subfigure}

\usepackage[authoryear]{natbib}

\usepackage[all,frame,ps,tips,line,poly,arrow,curve,import]{xy}
\SelectTips{cm}{10} \UseTips

\DeclareMathOperator*{\ext}{ext}

\DeclareMathOperator{\ad}{ad}
\DeclareMathOperator{\Ad}{Ad}
\DeclareMathOperator{\dexp}{dexp}
\DeclareMathOperator{\ddexp}{ddexp}

\def\vbar{\bar{\hat{v}}}

\newcommand{\bfi}{\bfseries\itshape}

\restylefloat{figure}

\newtheorem{theorem}{Theorem}[section]

\newtheorem{lemma}[theorem]{Lemma}

\newtheorem{remark}{Remark}[section]

\numberwithin{equation}{section}

\renewcommand{\paragraph}[1]{\vspace*{0.1in}\noindent\textbf{#1}}

\title{Generalized Galerkin Variational Integrators}
\author{Melvin Leok}
\address{Department of Mathematics, University of Michigan, Ann Arbor, MI 48109.}
\email{mleok@umich.edu}

\allowdisplaybreaks

\begin{document}

\begin{abstract}
We introduce generalized Galerkin variational integrators, which are a natural generalization of discrete variational mechanics, whereby the discrete action, as opposed to the discrete Lagrangian, is the fundamental object. This is achieved by approximating the action integral with appropriate choices of a finite-dimensional function space that approximate sections of the configuration bundle and numerical quadrature to approximate the integral. We discuss how this general framework allows us to recover higher-order Galerkin variational integrators, asynchronous variational integrators, and symplectic-energy-momentum integrators. In addition, we will consider function spaces that are not parameterized by field values evaluated at nodal points, which allows the construction of Lie group, multiscale, and pseudospectral variational integrators. The construction of pseudospectral variational integrators is illustrated by applying it to the (linear) Schr\"odinger equation. $G$-invariant discrete Lagrangians are constructed in the context of Lie group methods through the use of natural charts and interpolation at the level of the Lie algebra. The reduction of these $G$-invariant Lagrangians yield a higher-order analogue of discrete Euler--Poincar\'e reduction. By considering nonlinear approximation spaces, spatio-temporally adaptive variational integrators can be introduced as well.
\end{abstract}

\maketitle

\setcounter{tocdepth}{1} \tableofcontents

\section{Introduction}
We will review some of the previous work on discrete mechanics and their multisymplectic generalizations (see, for example,~\cite{MaPaSh1998, MaPeShWe2001}), before introducing a general formulation of discrete mechanics that recovers higher-order variational integrators (see, for example,~\cite{MaWe2001}), asynchronous variational integrators (see, for example,~\cite{LeMaOrWe2003}), as well symplectic-energy-momentum integrators (see, for example,~\cite{KaMaOr1999}).

While discrete variational integrators exhibit desirable properties such as symplecticity, momentum preservation, and good energy behavior, it does not address other important issues in numerical analysis, such as adaptivity and approximability. Generalized variational integrators are introduced with a view towards addressing such issues in the context of discrete variational mechanics.

By formulating the construction of a generalized variational integrator in terms of the choice of a finite-dimensional function space and a numerical quadrature scheme, we are able to draw upon the extensive literature on approximation theory and numerical quadrature to construct variational schemes that are appropriate for a larger class of problems. Within this framework, we will introduce multiscale, spatio-temporally adaptive, Lie group, and pseudospectral variational integrators.

\subsection{Standard Formulation of Discrete Mechanics}\index{discrete mechanics}
The standard formulation of discrete variational mechanics (see, for example, \cite{MaWe2001}) is to
consider the {\bfi discrete Hamilton's principle}\index{variational principle!discrete},
\[ \delta \mathbb{S}_d = 0,\]
where the {\bfi discrete action sum}, $\mathbb{S}_d:Q^{n+1}\rightarrow \mathbb{R}$, is given by
\[ \mathbb{S}_d(q_0,q_1,\ldots,q_n) = \sum_{i=0}^{n-1} L_d(q_i, q_{i+1}) .\]
The {\bfi discrete Lagrangian}\index{Lagrangian!discrete}, $L_d:Q\times Q\rightarrow \mathbb{R}$, is a generating function of the symplectic flow, and is an approximation to the {\bfi exact discrete Lagrangian}\index{Lagrangian!discrete!exact},
\[L_d^{\operatorname{exact}}(q_0,q_1)=\int_0^h L(q_{01}(t),\dot q_{01}(t)) dt,\]
where $q_{01}(0)=q_0,$ $q_{01}(h)=q_1,$ and $q_{01}$ satisfies the
Euler--Lagrange equation in the time interval $(0,h)$. The exact discrete Lagrangian is related to the Jacobi solution of the Hamilton--Jacobi equation. The discrete variational principle then yields the {\bfi discrete Euler--Lagrange (DEL)} equation,
\[ D_2 L_d(q_0,q_1)+D_1 L_d(q_1,q_2)=0,\]
which yields an implicit update map $(q_0,q_1)\mapsto(q_1,q_2)$ that is valid for initial conditions $(q_0,q_1)$ that as sufficiently close to the diagonal of $Q\times Q$.

\subsection{Multisymplectic Geometry}\index{multisymplectic!geometry}
The generalization of the variational principle to the setting of partial differential equations involves a multisymplectic formulation (see, for example, \cite{MaPaSh1998, MaPeShWe2001}). Here, the {\bfi base space} $\mathcal{X}$ consists of the independent variables, which are denoted by $(x^0,\ldots,x^n)$, where $x^0$ is time, and $x^1,\ldots,x^n$ are space variables.

The independent or field variables, denoted $(q^1,\ldots, q^m)$, form a fiber over each space-time basepoint. The set of independent variables, together with the field variables over them, form a fiber bundle, $\pi:Y\rightarrow\mathcal{X}$, referred to as the {\bfi configuration bundle}\index{multisymplectic!configuration bundle}. The configuration of the system is specified by giving the field values at each space-time point. More precisely, this can be represented as a section of $Y$ over $\mathcal{X}$, which is a continuous map $q:\mathcal{X}\rightarrow Y$, such that $\pi\circ q=1_{\mathcal{X}}$. This means that for every $x\in\mathcal{X}$, $q(x)$ is in the fiber over $x$, which is $\pi^{-1}(x)$.

In the case of ordinary differential equations, the Lagrangian is dependent on the position variable, and its time derivative, and the action integral is obtained by integrating the Lagrangian in time. In the multisymplectic case, the Lagrangian density is dependent on the field variables, and the derivatives of the field variables with respect to the space-time variables, and the action integral is obtained by integrating the Lagrangian density over a region of space-time.

The analogue of the tangent bundle $TQ$ in the multisymplectic setting is referred to as the {\bfi first jet bundle}\index{multisymplectic!first jet bundle} $J^1 Y$, which consists of the configuration bundle $Y$, together with the first derivatives of the field variables with respect to independent variables. We denote these as,
\[ {v^i}_j={q^i}_{,j}= \frac{\partial q^i}{\partial x^j},\]
for $i=1,\ldots, m$, and $j=0,\ldots, n$.

We can think of $J^1 Y$ as a fiber bundle over $\mathcal{X}$. Given a section $q:\mathcal{X}\rightarrow Y$, we obtain its {\bfi first jet extension}\index{multisymplectic!first jet extension}, $j^1 q:\mathcal{X}\rightarrow J^1 Y$, that is given by
\[ j^1 q (x^0,\ldots, x^n) = \left( x^0,\ldots, x^n, q^1,\ldots, q^m, {q^1}_{,1},\ldots, {q^m}_{,n} \right),\]
which is a section of the fiber bundle $J^1 Y$ over $\mathcal{X}$.

The Lagrangian density is a map $L:J^1 Y\rightarrow \Omega^{n+1}(\mathcal{X})$, and the action integral is given by
\[ \mathcal{S}(q) = \int_{\mathcal{X}} L(j^1 q), \]
and then Hamilton's principle states that
\[ \delta \mathcal{S} = 0.\]
We will see in the next subsection how this allows us to construct multisymplectic variational integrators.

\subsection{Multisymplectic Variational Integrator}\label{gvi:subsec:multisymplectic_integrator}\index{multisymplectic!variational integrator|see{variational integrator, multisymplectic}}\index{variational integrator!multisymplectic}
We introduce a multisymplectic variational integrator through the use of a simple but illustrative example. We consider a tensor product discretization of $(1+1)$-space-time, given by
\[
\begin{xy} 0;<10mm,0cm>:<0cm,10mm>::
(1,0)="a01", (3,0)="a02", (5,0)="a03", (0,1)="b01", (0,3)="b02",
(0,5)="b03", (1,6)="a11", (3,6)="a12", (5,6)="a13", (6,1)="b11",
(6,3)="b12", (6,5)="b13", (1,1)*@{*}, (1,3)*@{*}, (1,5)*@{*},
(3,1)*@{*}, (3,3)*@{*}, (3,5)*@{*}, (5,1)*@{*}, (5,3)*@{*},
(5,5)*@{*}, (1,-0.25)*\txt{$x_{i-1}$}, (3,-0.25)*\txt{$x_i$},
(5,-0.25)*\txt{$x_{i+1}$}, (-0.3,1)*\txt{$t_{j-1}$},
(-0.3,3)*\txt{$t_j$}, (-0.3,5)*\txt{$t_{j+1}$},
(3.3,4)*\txt{$\Delta t$}, (4,3.25)*\txt{$\Delta x$}, \ar@{-}
"a01";"a11", \ar@{-} "a02";"a12", \ar@{-} "a03";"a13", \ar@{-}
"b01";"b11", \ar@{-} "b02";"b12", \ar@{-} "b03";"b13"
\end{xy}
\]
and tensor product shape functions given by
\[\varphi_{i,j}(x,t) =
\raisebox{-10mm}{\begin{xy}
0;<10mm,0cm>:<0cm,10mm>::
(1.5,0.25)*\txt{$x_i$},
(2.5,0.25)*\txt{$x_{i+1}$},
(0.25,1.5)*\txt{$1$},
\ar@{-} (1.5,0.4);(1.5,0.6),
\ar@{-} (2.5,0.4);(2.5,0.6),
\ar@{-} (0.4,1.5);(0.6,1.5),
\ar (0,0.5);(3,0.5), \ar(0.5,0);(0.5,2), \ar@{-} (0.5,0.5);(1.5,1.5), \ar@{-} (1.5,1.5);(2.5,0.5)
\end{xy}}\,\,
\otimes\,\,
\raisebox{-10mm}{\begin{xy}
0;<10mm,0cm>:<0cm,10mm>::
(1.5,0.25)*\txt{$t_j$},
(2.5,0.25)*\txt{$t_{j+1}$},
(0.25,1.5)*\txt{$1$},
\ar@{-} (1.5,0.4);(1.5,0.6),
\ar@{-} (2.5,0.4);(2.5,0.6),
\ar@{-} (0.4,1.5);(0.6,1.5),
\ar (0,0.5);(3,0.5), \ar(0.5,0);(0.5,2), \ar@{-} (0.5,0.5);(1.5,1.5), \ar@{-} (1.5,1.5);(2.5,0.5)
\end{xy}}
\]
We construct the discrete Lagrangian as follows,
\[
L_d(q_{i,j}, q_{i+1,j}, q_{i,j+1}, q_{i+1,j+1})=\int_{[x_i,x_{i+1}]}\int_{[t_j,t_{j+1}]}
L\Big (j^1\Big (\sum\nolimits_{a=i}^{i+1}\sum\nolimits_{b=j}^{j+1} q_{a,b}\varphi_{a,b}\Big )\Big),
\]
where $q_{i,j}\simeq q(i\Delta x, j\Delta t)$.

Consider a variation that varies the value of $q_{i,j}$, and leaves the other degrees of freedom fixed, then we obtain the following discrete Euler--Lagrange equation,
\begin{align*}
D_1 L_d(q_{i,j},q_{i+1,j},q_{i,j+1},q_{i+1,j+1}) +&D_2 L_d(q_{i-1,j},q_{i,j},q_{i-1,j+1},q_{i,j+1})\\
+&D_3 L_d(q_{i,j-1},q_{i+1,j-1},q_{i,j},q_{i+1,j})\\
+&D_4 L_d(q_{i-1,j-1},q_{i,j-1},q_{i-1,j},q_{i,j})=0\, ,
\end{align*}
In general, when we take variations in a degree of freedom, we will have a discrete Euler--Lagrange equation that involves all terms in the discrete action sum that are associated with regions in space-time that overlap with the support of the shape function associated with that degree of freedom.

\section{Generalized Galerkin Variational Integrators}\index{variational integrator!generalized}
There are a few essential observations that go into constructing a general framework that encompasses the prior work on variational integrators, asynchronous variational integrators, and symplectic-energy-momentum integrators, while yielding generalizations that allow the construction of multiscale, spatio-temporally adaptive, Lie group, and pseudospectral variational integrators.

The first is that a generalized Galerkin variational integrator involves the choice of a finite-dimensional function space that discretizes a section of the configuration bundle, and the second is that we approximate the action integral though a numerical quadrature scheme to yield a discrete action sum. The discrete variational equations we obtain from this discrete variational principle are simply the Karush--Kuhn--Tucker (KKT) conditions (see, for example,~\cite{NoWr1999}) with respect to the degrees of freedom that generate the finite-dimensional function space.

To recap, the choices which are made in discretizing a variational problem are:
\begin{enumerate}
\item A finite-dimensional function space to represent sections of the configuration bundle.\index{variational integrator!function space}
\item A numerical quadrature scheme to evaluate the action integral.\index{variational integrator!numerical quadrature}
\end{enumerate}
Given these two choices, we obtain an expression for the discrete action in terms of the degrees of freedom, and the KKT conditions for the discrete action to be stationary with respect to variations in the degrees of freedom yield the generalized discrete Euler--Lagrange equations.

Current variational integrators are based on piecewise
polynomial interpolation, with function spaces that are parameterized by
the value of the field variables at nodal points, and a set of
internal points. By relaxing the condition that the interpolation
is piecewise, we will be able to consider pseudospectral
discretizations, and by relaxing the condition that the
parameterization is in terms of field values, we will be able to
consider Lie group variational integrators. By considering
shape functions motivated by multiscale finite elements (see, for
example,~\cite{HoWu1999, EfHoWu2000, ChHo2003}), we will obtain
multiscale variational integrators. And by generalizing the
approach used in symplectic-energy-momentum integrators, and
considering nonlinear approximation spaces (see, for
example,~\cite{DeVore1998}), we will be able to introduce
spatio-temporally adaptive variational integrators.

As we will see, there is nothing canonical about the form of the discrete Euler--Lagrange equations, or the notion that the discrete Lagrangian is a map $L_d:Q\times Q\rightarrow \mathbb{R}$. These expressions arise because of the finite-dimensional function space which has been chosen.

\subsection{Special Cases of Generalized Galerkin Variational Integrators}
In this subsection, we will show how higher-order Galerkin variational integrators, multisymplectic variational integrators, and symplectic-energy-momentum integrators are all special cases of generalized Galerkin variational integrators.

\paragraph{Higher-Order Galerkin Variational Integrators.}\index{variational integrator!higher-order} In the case of higher-order Galerkin variational integrators, we have chosen a piecewise interpolation for each time interval $[0,h]$, that is parameterized by control points $q_0^1,\ldots q_0^s$, corresponding to the value of the curve at a set of control times $0=d_0<d_1<\ldots<d_{s-1}<d_s=1$. The interpolation within each interval $[0,h]$ is given by the unique degree $s$ polynomial $q_d(t;q_0^\nu,h)$, such that $q_d(d_\nu h) = g_0^\nu$, for $\nu=0,\ldots, s$.

\begin{figure}[htbp]
\[
\begin{xy} 0;<30mm,0cm>:<0cm,10mm>::
(-0.1,4.6)*\txt{$Q$},
(2.65,-0.3)*\txt{$t$},
(-0.1,-0.3)*\txt{$0$},
(0.5,-0.3)*\txt{$d_1 h$},
(1,-0.3)*\txt{$d_2 h$},
(1.5,-0.3)*\txt{$d_{s-2} h$},
(2,-0.3)*\txt{$d_{s-1} h$},
(2.5,-0.3)*\txt{$h$},
(0,2)*@{*};
(1,1)*@{*}
**\crv{(0.5,3)}
?(.5)*@{*} * !LD!/^-5pt/{q_0^1},
(1,1);(1.5,2) **\crv{~*=<4pt>{.} (1.2,0)},
(1.5,2)*@{*};(2.5,4)*@{*}
**\crv{(2,4.3)}
?(.5)*@{*} *!LD!/^-5pt/{q_0^{s-1}},
(0.1,1.7)* !LD!{q_0^0},
(1,1.3)* !LD!{q_0^2},
(1.45,2.2)* !LD!{q_0^{s-2}},
(2.5,4.2)* !LD!{q_0^s},
\ar@{-} (-0.1,0);(1,0),
\ar@{.} (1,0);(1.5,0),
\ar (1.5,0);(2.6,0),
\ar (0,-0.3);(0,4.5),
\ar@{-} (0.5,-0.1);(0.5,0.1),
\ar@{-} (1,-0.1);(1,0.1),
\ar@{-} (1.5,-0.1);(1.5,0.1),
\ar@{-} (2,-0.1);(2,0.1),
\ar@{-} (2.5,-0.1);(2.5,0.1),
\end{xy}
\]
\caption{Polynomial interpolation used in higher-order Galerkin variational integrators.}
\end{figure}
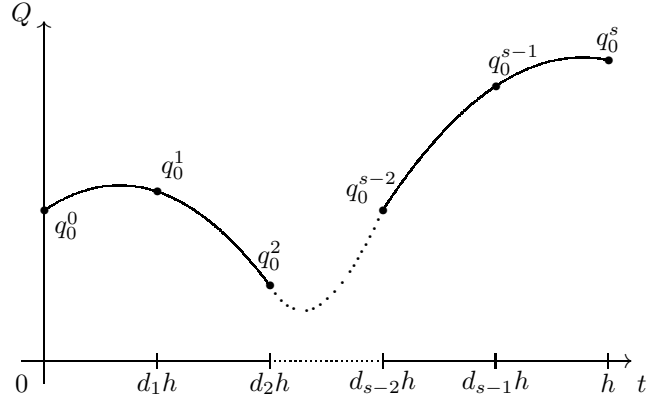

By an appropriate choice of quadrature scheme, we can break up the action integral into pieces, which we denote by
\[ \mathbb{S}_d^i(q_i^\nu)\approx \int_0^h L(q_d(t;q_i^\nu,h),\dot q_d(t;q_i^\nu,h)) dt. \]
If we further require that the piecewise defined curve is continuous at the node points, we obtain the following augmented discrete action,
\begin{align*}
\mathbb{S}_d\left(\{q_i^\nu\}_{\stackrel{\scriptstyle i=0,\ldots, N-1}{\nu=0,\ldots,s}}\right)
&=\sum_{i=0}^{N-1} \mathbb{S}_d^i\left(\{q_i^\nu\}_{\nu=0}^s\right) - \sum_{i=0}^{N-2}\lambda_i (q_i^s-q_{i+1}^0).
\end{align*}
The discrete action is stationary when
\begin{align*}
\frac{\partial \mathbb{S}_d^i}{\partial q_i^s}(q_i^\nu) & = \lambda_i,\\
\frac{\partial \mathbb{S}_d^i}{\partial q_{i+1}^0}(q_{i+1}^\nu) & = -\lambda_i,\\
\frac{\partial \mathbb{S}_d^i}{\partial q_i^j}(q_i^\nu) & = 0.
\end{align*}
We can identify the pieces of the discrete action with the discrete Lagrangian, $L_d:Q\times Q\rightarrow \mathbb{R}$, by setting
\begin{equation}
L_d(q_i, q_{i+1}) = \mathbb{S}_d^i(q_i^\nu),\label{gvi:eqn:discrete_lagrangian_discrete_action}
\end{equation}
where, $q_i^0=q_i$, $q_i^s=q_{i+1}$, and the other $q_i^j$'s are defined implicitly by the system of equations
\begin{equation}
\frac{\partial \mathbb{S}_d^i}{\partial q_i^j}(q_i^\nu) = 0,\label{gvi:eqn:high_order_Ld_stationary}
\end{equation}
for $\nu=1,\ldots, s-1$.
Once we have made this identification, we have that
\begin{align*}
-D_1 L_d(q_{i+1},q_{i+2})
&= -\frac{\mathbb{S}_d^{i+1}}{\partial q_{i+1}^0}(q_{i+1}^\nu)\\
&=\lambda_i \\
&= \frac{\partial \mathbb{S}_d^i}{\partial q_i^s}(q_i^\nu)\\
&= D_2 L_d(q_i, q_{i+1}),
\end{align*}
from which we recover the discrete Euler--Lagrange equation,
\[ D_1 L_d(q_{i+1},q_{i+2}) + D_2 L_d(q_i, q_{i+1}) =0. \]
The DEL equations, together with the definition of the discrete Lagrangian, given in Equations~\ref{gvi:eqn:discrete_lagrangian_discrete_action} and \ref{gvi:eqn:high_order_Ld_stationary}, yield a higher-order Galerkin variational integrator. As we have shown, it is only because the interpolation was piecewise that we were able to decompose the equations into a DEL equation, and a set of equations that define the discrete Lagrangian in terms of conditions on the internal control points.

\paragraph{Multisymplectic Variational Integrators.}\index{variational integrator!multisymplectic}
Recall the example of a multisymplectic variational integrator we introduced in \S\ref{gvi:subsec:multisymplectic_integrator}, where we used tensor product linear shape functions with localized supports to discretize the configuration bundle. The discrete action is then given by
\begin{align*}
\mathbb{S}_d\left(\{q_{a,b}\}_{\stackrel{\scriptstyle a=0,\ldots,M-1}{b=0,\ldots, N-1}}\right)
&=\int_{[x_0,x_M]}\int_{[t_0,t_N]} L \Big(j^1\Big(\sum_{a=0}^N \sum_{b=0}^N q_{a,b} \varphi_{a,b}\Big)\Big)\\
&=\sum_{i=0}^{M-1} \sum_{j=0}^{N-1}\int_{[x_i,x_{i+1}]}\int_{[t_j,t_{j+1}]} L \Big(j^1\Big(\sum_{a=0}^N \sum_{b=0}^N q_{a,b} \varphi_{a,b}\Big)\Big)\\
&=\sum_{i=0}^{M-1} \sum_{j=0}^{N-1}\int_{[x_i,x_{i+1}]}\int_{[t_j,t_{j+1}]} L \Big(j^1\Big(\sum_{a=i}^{i+1} \sum_{b=j}^{j+1} q_{a,b} \varphi_{a,b}\Big)\Big)\, ,
\end{align*}
where we first decomposed the integral into pieces, and then used the local support of the shape functions to simplify the inner sums over $q_{a,b}\varphi_{a,b}$. As before, it is due to the local support of the shape functions that we can express the discrete action as
\[ \mathbb{S}_d = \sum_{i=0}^{M-1} \sum_{j=0}^{N-1} L_d (q_{i,j}, q_{i+1,j}, q_{i,j+1}, q_{i+1,j+1}),
\]
where
\[
L_d(q_{i,j}, q_{i+1,j}, q_{i,j+1}, q_{i+1,j+1})=\int_{[x_i,x_{i+1}]}\int_{[t_j,t_{j+1}]}
L\Big (j^1\Big (\sum\nolimits_{a=i}^{i+1}\sum\nolimits_{b=j}^{j+1} q_{a,b}\varphi_{a,b}\Big )\Big),
\]
This localized support is also the reason why the discrete Euler--Lagrange equation in this case consists of four terms, since to each degree of freedom, there are four other degrees of freedom that have shape functions with overlapping support. In the case of ordinary differential equations, this was two. In general, for tensor product meshes of $(n+1)$-space-times, with tent function shape functions, the number of terms in the discrete Euler--Lagrange equation will be $2^{n+1}$. In contrast, for pseudospectral variational integrators with $m$ spatial degrees of freedom per time level, and $k$ degrees of freedom per piecewise polynomial in time, each of the $m(k-1)$ discrete Euler--Lagrange equations will involve $m(k-1)$ terms.

This is simply a reflection of the fact that shape functions with compact support yield schemes with banded matrix structure, whereas pseudospectral and spectral methods tend to yield fuller matrices. The payoff in using pseudospectral and spectral methods for problems with smooth or analytic solutions is due to the approximation theoretic property that these solutions are approximated at an exponential rate of accuracy by spectral expansions.

\paragraph{Symplectic-Energy-Momentum Integrators.}\index{variational integrator!symplectic-energy-momentum}
In the case of symplectic-energy-momentum integrators, the degrees
of freedom involve both the base variables and the field
variables. We will first derive a second-order symplectic-energy
momentum integrator, and in the next section, we will derive a
higher-order generalization. We choose a piecewise linear
interpolation for our configuration bundle. Each piece is
parameterized by the endpoint values of the field variable
$q_i^0$, $q_i^1$, and the endpoint times $h_i^0$, $h_i^1$, and we
approximate the action integral using the midpoint rule. To ensure
continuity, we require that, $q_i^1=q_{i+1}^0$, and
$h_i^1=h_{i+1}^0$.

\begin{remark}The approach of allowing each piecewise defined curve to
float around freely, and imposing the continuity conditions using
Lagrange multipliers was used in \cite{LaWe2003} to unify the formulation of
discrete variational mechanics and optimal control. Through the use of a primal-dual formalism, discrete analogues of Hamiltonian mechanics and the Hamilton--Jacobi equation were also introduced.
\end{remark}

This yields the following discrete action,
\[ \mathbb{S}_d = \sum_{i=0}^{N-1}(h_i^1-h_i^0) L\left(\frac{q_i^0+q_i^1}{2},\frac{q_i^1-q_i^0}{h_i^1-h_1^0}\right)-\sum_{i=0}^{N-2} \lambda_i (q_i^1-q_{i+1}^0) - \sum_{i=0}^{N-2}
\omega_i(h_i^1-h_{i+1}^0).
\]
To simplify the expressions, we define
\begin{align*}
h_i&\equiv h_i^1-h_i^0,\\
q_{i+\frac{1}{2}}&\equiv\frac{q_i^0+q_i^1}{2},\\
\dot q_{i+\frac{1}{2}}&\equiv\frac{q_i^1-q_i^0}{h_i}.
\end{align*}
Then, the variational equations are given by
\begin{align*}
0 =& L\left(q_{i+\frac{1}{2}},\dot q_{i+\frac{1}{2}}\right) -
h_i\frac{\partial L}{\partial \dot q}\left(q_{i+\frac{1}{2}},\dot
q_{i+\frac{1}{2}}\right)\frac{1}{h_i}\dot
q_{i+\frac{1}{2}}-\omega_i, &\text{for }i=0,\ldots, N-2,\\
0 =&  -L\left(q_{i+\frac{1}{2}},\dot q_{i+\frac{1}{2}}\right) +
h_i\frac{\partial L}{\partial \dot q}\left(q_{i+\frac{1}{2}},\dot
q_{i+\frac{1}{2}}\right)\frac{1}{h_i}\dot
q_{i+\frac{1}{2}}+\omega_{i-1}, &\text{for }i=1,\ldots, N-1,\\
0 =& h_i\left[\frac{\partial L}{\partial
q}\left(q_{i+\frac{1}{2}},\dot
q_{i+\frac{1}{2}}\right)\frac{1}{2}+\frac{\partial L}{\partial
\dot q}\left(q_{i+\frac{1}{2}},\dot
q_{i+\frac{1}{2}}\right)\frac{1}{h_i}\right]-\lambda_i, &\text{for
}i=0,\ldots,N-2,\\
0 =& h_i\left[\frac{\partial L}{\partial
q}\left(q_{i+\frac{1}{2}},\dot
q_{i+\frac{1}{2}}\right)\frac{1}{2}-\frac{\partial L}{\partial
\dot q}\left(q_{i+\frac{1}{2}},\dot
q_{i+\frac{1}{2}}\right)\frac{1}{h_i}\right]+\lambda_{i-1},
&\text{for
}i=1,\ldots,N-1,\\
q_i^1 =& q_{i+1}^0, &\text{for } i=0,\ldots,N-2,\\
h_i^1 =& h_{i+1}^0, &\text{for } i=0,\ldots,N-2.
\end{align*}
If we define the discrete Lagrangian to be
\[ L_d(q_i,q_{i+1},h_i) \equiv h_i
L\left(\frac{q_i+q_{i+1}}{2},\frac{q_{i+1}-q_i}{h_i}\right),\]
and the discrete energy to be
\begin{align*}
E_d(q_i,q_{i+1},h_i)
&\equiv -D_3 L_d(q_i,q_{i+1},h_i)\\
&=-L\left(\frac{q_i+q_{i+1}}{2},\frac{q_{i+1}-q_i}{h_i}\right)+
h_i\frac{\partial L}{\partial \dot
q}\left(\frac{q_i+q_{i+1}}{2},\frac{q_{i+1}-q_i}{h_i}\right)\frac{1}{h_i}\frac{q_{i+1}-q_i}{h_i},
\end{align*}
and identify points as follows,
\begin{align*}
q_i &= q_{i-1}^1 = q_i^0,\\
h_i &= h_{i-1}^1 = h_i^0,
\end{align*}
we obtain
\begin{align*}
0 =& E_d(q_i,q_{i+1},h_i)-\omega_i, &\text{for }i=0,\ldots, N-2,\\
0 =& -E_d(q_i,q_{i+1},h_i)+\omega_{i-1}, &\text{for }i=1,\ldots, N-1,\\
0 =& D_2 L_d(q_i,q_{i+1},h_i)-\lambda_i, &\text{for
}i=0,\ldots,N-2,\\
0 =& D_1 L_d(q_i,q_{i+1},h_i)+\lambda_{i-1}, &\text{for
}i=1,\ldots,N-1.
\end{align*}
After eliminating the Lagrange multipliers, we obtain the
conservation of discrete energy equation,
\[
E_d(q_i, q_{i+1},h_i)= E_d(q_{i+1},q_{i+2},h_{i+1}),
\]
and the discrete Euler--Lagrange equation,
\[
D_2 L_d (q_i,q_{i+1},h_i) + D_1 L_d (q_{i+1},q_{i+2},h_{i+1}) = 0.
\]
This recovers the results obtained in \cite{KaMaOr1999}, but the
derivation is new. In \S\ref{gvi:sec:higher_order_sym_eng_mom}, we will see how this
derivation can be generalized to yield a higher-order scheme.

\paragraph{Discrete Action as the Fundamental Object.}\index{variational integrator!discrete action} The message of this section is that the discrete action is the fundamental object in discrete mechanics, as opposed to the discrete Lagrangian. In instances whereby the shape function associated with individual degrees of freedom have localized supports, it is possible to decompose the discrete action into terms that can be identified with discrete Lagrangians. While this approach might seem artificial at first, we will find that in discussing pseudospectral variational integrators, it does not make sense to break up the discrete action into individual pieces.

\section{Lie Group Variational Integrators}\index{variational integrator!Lie group}
In this section, we will introduce higher-order Lie group variational integrators. The basic idea behind all Lie group techniques is to express the update map of the numerical scheme in terms of the exponential map,
\[ g_1 = g_0 \exp(\xi_{01})\, ,\]
and thereby reduce the problem to finding an appropriate Lie
algebra element $\xi_{01}\in\mathfrak{g}$, such that the update
scheme has the desired order of accuracy. This is a desirable
reduction, as the Lie algebra is a vector space, and as such the
interpolation of elements can be easily defined. In our
construction, the interpolatory method we use on the Lie group
relies on interpolation at the level of the Lie algebra.

For a more in depth review of Lie group methods, please refer to~\cite{IsMuNoZa2000}.
In the case of variational Lie group methods, we will express the variational problem in terms of finding Lie algebra elements, such that the discrete action is stationary.

As we will consider the reduction of these higher-order Lie group integrators in the next section, we will chose a construction that yields a $G$-invariant discrete Lagrangian whenever the continuous Lagrangian is $G$-invariant. This is achieved through the use of $G$-equivariant interpolatory functions, and in particular, natural charts on $G$.

\subsection{Galerkin Variational Integrators}\index{variational integrator!Galerkin}
We first recall the construction of higher-order Galerkin
variational integrators, as originally described in \cite{MaWe2001}.
Given a Lie group $G$, the associated {\bfi state space} is
given by the tangent bundle $TG$. In addition, the dynamics on $G$
is described by a {\bfi Lagrangian}, $L:TG\rightarrow\mathbb{R}$.
Given a time interval $[0,h]$, the {\bfi path space} is defined
to be
\[\mathcal{C}(G)=\mathcal{C}([0,h],G)=\{g:[0,h]\rightarrow G\mid g \text{ is a } C^2
\text{ curve}\},\] and the {\bfi action map},
$\mathfrak{S}:\mathcal{C}(G)\rightarrow\mathbb{R}$, is given by
\[\mathfrak{S}(g)\equiv \int_0^h L(g(t),\dot g(t)) dt.\]

We approximate the action map, by numerical quadrature, to yield $\mathfrak{S}^s:\mathcal{C}([0,h],G)\rightarrow
\mathbb{R}$,
\[\mathfrak{S}^s(g)\equiv h \sum_{i=1}^s b_i L(g(c_i h), \dot g(c_i
h)),\] where $c_i\in[0,1]$, $i=1,\ldots, s$ are the quadrature
points, and $b_i$ are the quadrature weights.

Recall that the discrete Lagrangian should be an approximation of the form
\[L_d(g_0,g_1,h)\approx \ext_{g\in\mathcal{C}([0,h],G),\\g(0)=g_0, g(h)=g_1} \mathfrak{S}(g)\, .\]
If we restrict the extremization procedure to the subspace spanned
by the interpolatory function that is parameterized by $s+1$
internal points, $\varphi:G^{s+1}\rightarrow
\mathcal{C}([0,h],G)$, we obtain the following discrete
Lagrangian,
\begin{align*}
L_d(g_0,g_1) &= \ext_{g^\nu\in G; g^0=g_0; g^s=g_1} \mathfrak{S}
(T\varphi(g^\nu;\cdot))\\
&= \ext_{g^\nu\in G; g^0=g_0; g^s=g_1} h\sum_{i=1}^s b_i
L(T\varphi(g^\nu;c_i h)).
\end{align*}
The interpolatory function is $G$-equivariant if
\[\varphi(g g^\nu;t)=g \varphi( g^\nu;t ).\]
\begin{lemma}\label{gvi:lemma:invariant_Ld}
If the interpolatory function $\varphi(g^\nu;t)$ is
$G$-equivariant, and the Lagrangian, $L:TG\rightarrow\mathbb{R}$, is
$G$-invariant, then the discrete Lagrangian, $L_d:G\times
G\rightarrow \mathbb{R}$, given by
\[
L_d(g_0,g_1)= \ext_{g^\nu\in G; g^0=g_0; g^s=g_1} h\sum_{i=1}^s
b_i L(T\varphi(g^\nu;c_i h)),
\]
is $G$-invariant.
\end{lemma}
\begin{proof}
\begin{align*}
L_d(g g_0, g g_1) &= \ext_{\tilde g^\nu\in G; \tilde g^0=g g_0;
\tilde g^s=g g_1}
h\sum_{i=1}^s b_i L(T\varphi(\tilde g^\nu;c_i h)),\\
&= \ext_{g^\nu\in g^{-1}G; g^0=g_0; g^s=g_1}
h\sum_{i=1}^s b_i L(T\varphi( g g^\nu;c_i h)),\\
&= \ext_{g^\nu\in G; g^0=g_0; g^s=g_1}
h\sum_{i=1}^s b_i L(TL_g\cdot T\varphi(g^\nu;c_i h)),\\
&= \ext_{g^\nu\in G; g^0=g_0; g^s=g_1}
h\sum_{i=1}^s b_i L(T\varphi(g^\nu;c_i h)),\\
&= L_d(g_0,g_1),
\end{align*}
where we used the $G$-equivariance of the interpolatory function
in the third equality, and the $G$-invariance of the Lagrangian in
the forth equality.
\end{proof}

\begin{remark}
While $G$-equivariant interpolatory functions provide a computationally efficient method of constructing $G$-invariant discrete Lagrangians, we can construct a $G$-invariant discrete Lagrangian (when $G$ is compact) by averaging an arbitrary discrete Lagrangian. In particular, given a discrete Lagrangian $L_d:Q\times Q\rightarrow \mathbb{R}$, the averaged discrete Lagrangian, given by
\[ \bar{L}_d(q_0,q_1) = \frac{1}{|G|}\int_{g\in G} L_d(g q_0, g q_1) dg\]
is $G$-equivariant. Therefore, in the case of compact symmetry groups, a $G$-invariant discrete Lagrangian always exists.
\end{remark}

\subsection{Natural Charts}\index{natural charts}\label{gvi:subsec:natural_charts}
Following the construction in \cite{MaPeSh1999}, we use the group
exponential map at the identity, $\exp_e:\mathfrak{g}\rightarrow
G$, to construct a $G$-equivariant interpolatory function, and a
higher-order discrete Lagrangian. As shown in
Lemma~\ref{gvi:lemma:invariant_Ld}, this construction yields a
$G$-invariant discrete Lagrangian if the Lagrangian itself is
$G$-invariant.

In a finite-dimensional Lie group $G$, $\exp_e$ is a local
diffeomorphism, and thus there is an open neighborhood $U\subset
G$ of $e$ such that $\exp_e^{-1}:U\rightarrow
\mathfrak{u}\subset\mathfrak{g}$. When the group acts on the left,
we obtain a chart $\psi_g:L_g U \rightarrow \mathfrak{u}$ at $g\in
G$ by
\[\psi_g=\exp_e^{-1}\circ L_{g^{-1}}.\]

\begin{lemma}
The interpolatory function given by
\[\varphi(g^\nu;\tau h)=\psi_{g^0}^{-1}
\Big(\sum\nolimits_{\nu=0}^s
\psi_{g^0}(g^\nu)\tilde{l}_{\nu,s}(\tau)\Big),\] is
$G$-equivariant.
\end{lemma}
\begin{proof}
\begin{align*}
\varphi(gg^\nu;\tau h) &= \psi_{(gg^0)}^{-1}
\Big(\sum\nolimits_{\nu=0}^s
\psi_{gg^0}(gg^\nu)\tilde{l}_{\nu,s}(\tau)\Big)\\
&= L_{gg^0} \exp_e \Big(\sum\nolimits_{\nu=0}^s
\exp_e^{-1}((gg^0)^{-1}(gg^\nu))\tilde{l}_{\nu,s}(\tau)\Big)\\
&= L_g L_{g^0} \exp_e \Big(\sum\nolimits_{\nu=0}^s
\exp_e^{-1}((g^0)^{-1}g^{-1}gg^\nu)\tilde{l}_{\nu,s}(\tau)\Big)\\
&= L_g \psi_{g^0}^{-1} \Big(\sum\nolimits_{\nu=0}^s
\exp_e^{-1}\circ L_{(g^0)^{-1}}(g^\nu)\tilde{l}_{\nu,s}(\tau)\Big)\\
&= L_g \psi_{g^0}^{-1} \Big(\sum\nolimits_{\nu=0}^s
\psi_{g^0}(g^\nu)\tilde{l}_{\nu,s}(\tau)\Big)\\
&= L_g \varphi(g^\nu;\tau h).\qedhere
\end{align*}
\end{proof}
\begin{remark}
In the proof that $\varphi$ is $G$-equivariant, it was important
that the base point for the chart should transform in the same way
as the internal points $g^\nu$. As such, the interpolatory
function will be $G$-equivariant for a chart that it based at any
one of the internal points $g^\nu$ that parameterize the function,
but will not be $G$-equivariant if the chart is based at a fixed
$g\in G$. Without loss of generality, we will consider the case
when the chart is based at the first point $g_0$.
\end{remark}

We will now consider a discrete Lagrangian based on the use of interpolation in a natural chart, which is given by
\[ L_d(g_0, g_1) = \ext_{g^\nu\in G; g^0=g_0; g^s = g_0^{-1}g_1} h \sum_{i=1}^s b_i L(T\varphi (\{g^\nu\}_{\nu=0}^s;c_i h)\, .\]
To further simplify the expression, we will express the extremal in terms of the Lie algebra elements $\xi^\nu$ associated with the $\nu$-th control point. This relation is given by
\[\xi^\nu=\psi_{g_0}(g^\nu)\, ,\]
and the interpolated curve in the algebra is given by
\[\xi(\xi^\nu;\tau h)=\sum_{\kappa=0}^s \xi^\kappa \tilde l_{\kappa,s}(\tau),\]
which is related to the curve in the group,
\[g(g^\nu;\tau h)=g_0 \exp(\xi(\psi_{g_0}(g^\nu);\tau h)).\]
The velocity $\dot \xi=g^{-1} \dot g$ is given by
\[\dot\xi(\tau h)=g^{-1}\dot g(\tau h)=\frac{1}{h}\sum_{\kappa=0}^s \xi^\kappa \dot{\tilde
l}_{\kappa,s}(\tau).\]
Using the standard formula for the
derivative of the exponential,
\[T_\xi \exp= T_e L_{\exp(\xi)}\cdot \dexp_{\ad_{\xi}},\]
where
\[\dexp_w=\sum_{n=0}^\infty \frac{w^n}{(n+1)!},\]
we obtain the following expression for discrete Lagrangian,
\begin{align*}
L_d(g_0, g_1) &= \ext_{\xi^\nu\in\mathfrak{g};\xi^0=0;\xi^s=\psi_{g_0}(g_1)}
h \sum_{i=1}^s b_i L\Big(L_{g_0} \exp(\xi(c_i h)),\\
&\hspace*{2in} T_{\exp(\xi(c_i h))}L_{g_0}\cdot T_e
L_{\exp(\xi(c_i h))}\cdot \dexp_{\ad_{\xi(c_i h)}} (\dot \xi(c_i
h))\Big)\, .
\end{align*}
More explicitly, we can compute the conditions on the Lie algebra elements for the expression above to be extremal. This implies that
\[ L_d(g_0, g_1) = h \sum_{i=1}^s b_i L\Big(L_{g_0} \exp(\xi(c_i h)), T_{\exp(\xi(c_i h))}L_{g_0}\cdot T_e L_{\exp(\xi(c_i h))}\cdot \dexp_{\ad_{\xi(c_i h)}} (\dot \xi(c_i h))\Big)
\]
with $\xi^0=0$, $\xi^s=\psi_{g_0}(g_1)$, and the other Lie algebra elements implicitly defined by
\begin{align*}
0 &= h\sum_{i=1}^s b_i \left [\frac{\partial L}{\partial g} (c_i h) T_{\exp(\xi(c_i h))} L_{g_0}\cdot T_e L_{\exp(\xi(c_i h))} \cdot \dexp_{\ad_{\xi(c_i h)}} \tilde{l}_{\nu,s}(c_i) \right.\\
&\qquad\qquad\qquad+ \left . \frac{1}{h}\frac{\partial L}{\partial
\dot g}(c_i h) T^2_{\exp(\xi(c_i h))} L_{\exp(\xi(c_i h))}\cdot
T^2_e L_{\exp(\xi(c_i h))}\cdot \ddexp_{\ad_{\xi(c_i
h)}}\dot{\tilde{l}}_{\nu,s}(c_i)\right ]\, ,
\end{align*}
for $\nu=1,\ldots,s-1$, and where
\[ \ddexp_w = \sum_{n=0}^\infty \frac{w^n}{(n+2)!}\, . \]
This expression for the higher-order discrete Lagrangian, together with the discrete Euler--Lagrange equation,
\[ D_2 L_d(g_0,g_1) + D_1 L_d(g_1, g_2)=0\, ,\]
yields a {\bfi higher-order Lie group variational integrator}.

\section{Higher-Order Discrete Euler--Poincar\'e Equations}\index{variational integrator!Euler--Poincar\'e}\index{reduction!Euler--Poincar\'e}
In this section, we will apply discrete Euler--Poincar\'e reduction (see, for example,~\cite{MaPeSh1999}) to the Lie group variational integrator we derived previously, to construct a higher-order generalization of discrete Euler--Poincar\'e reduction.

\subsection{Reduced Discrete Lagrangian}
We first proceed by computing an expression for the reduced
discrete Lagrangian in the case when the Lagrangian is
$G$-invariant. Recall that our discrete Lagrangian uses
$G$-equivariant interpolation, which, when combined with the
$G$-invariance of the Lagrangian, implies that the discrete
Lagrangian is $G$-invariant as well. We compute the reduced
discrete Lagrangian,
\begin{align*}
l_d(g_0^{-1}g_1)
&\equiv L_d(g_0,g_1)\\
&= L_d(e, g_0^{-1}g_1)\\
&= \ext_{\xi^\nu\in\mathfrak{g};\xi^0=0;\xi^s=\log(g_0^{-1}g_1)}
h \sum_{i=1}^s b_i L\Big(L_{e} \exp(\xi(c_i h)),\\
&\hspace*{2.1in} T_{\exp(\xi(c_i h))}L_{e}\cdot T_e L_{\exp(\xi(c_i h))}\cdot \dexp_{\ad_{\xi(c_i h)}} (\dot \xi(c_i h))\Big)\\
&= \ext_{\xi^\nu\in\mathfrak{g};\xi^0=0;\xi^s=\log(g_0^{-1}g_1)}
h\sum_{i=1}^s b_i L\Big(\exp(\xi(c_i h)), T_e L_{\exp(\xi(c_i
h))}\cdot \dexp_{\ad_{\xi(c_i h)}}(\dot\xi(c_i h))\Big)\, .
\end{align*}
Setting $\xi^0=0$, and $\xi^s=\log(g_0^{-1}g_1)$, we can solve the
stationarity conditions for the other Lie algebra elements
$\{\xi^\nu\}_{\nu=1}^{s-1}$ using the following implicit system of
equations,
\begin{align*}
0 &= h\sum_{i=1}^s b_i \left [\frac{\partial L}{\partial g} (c_i h) T_e L_{\exp(\xi(c_i h))} \cdot \dexp_{\ad_{\xi(c_i h)}} \tilde{l}_{\nu,s}(c_i) \right.\\
&\qquad\qquad\qquad+ \left . \frac{1}{h}\frac{\partial L}{\partial
\dot g}(c_i h) T^2_e L_{\exp(\xi(c_i h))}\cdot
\ddexp_{\ad_{\xi(c_i h)}}\dot{\tilde{l}}_{\nu,s}(c_i)\right ]
\end{align*}
where $\nu=1,\ldots, s-1$.

This expression for the reduced discrete Lagrangian is not fully
satisfactory however, since it involves the Lagrangian, as opposed
to the reduced Lagrangian. If we revisit the expression for the
reduced discrete Lagrangian,
\[
l_d(g_0^{-1}g_1) =
\ext_{\xi^\nu\in\mathfrak{g};\xi^0=0;\xi^s=\log(g_0^{-1}g_1)}
h\sum_{i=1}^s b_i L\Big(\exp(\xi(c_i h)), T_e L_{\exp(\xi(c_i
h))}\cdot \dexp_{\ad_{\xi(c_i h)}}(\dot\xi(c_i h))\Big)\, ,\] we
find that by $G$-invariance of the Lagrangian, each of the terms
in the summation,
\[ L\Big(\exp(\xi(c_i h)), T_e L_{\exp(\xi(c_i
h))}\cdot \dexp_{\ad_{\xi(c_i h)}}(\dot\xi(c_i h))\Big)\, ,\] can
be replaced by
\[ l\Big(\dexp_{\ad_{\xi(c_i h)}}(\dot\xi(c_i h))\Big)\, ,\]
where $l:\mathfrak{g}\rightarrow \mathbb{R}$ is the {\bfi reduced
Lagrangian} given by
\[ l(\eta)= L(L_{g^{-1}}g,TL_{g^{-1}}\dot g) = L(e,\eta), \]
where $\eta=TL_{g^{-1}}\dot g\in\mathfrak{g}$.

From this observation, we have an expression for the reduced
discrete Lagrangian in terms of the reduced Lagrangian,
\[l_d(g_0^{-1}g_1) =
\ext_{\xi^\nu\in\mathfrak{g};\xi^0=0;\xi^s=\log(g_0^{-1}g_1)}
h\sum_{i=1}^s b_i l\Big(\dexp_{\ad_{\xi(c_i h)}}(\dot\xi(c_i
h))\Big)\, .\] As before, we set $\xi^0=0$, and
$\xi^s=\log(g_0^{-1}g_1)$, and solve the stationarity conditions
for the other Lie algebra elements $\{\xi^\nu\}_{\nu=1}^{s-1}$
using the following implicit system of equations,
\begin{align*}
0 &= h\sum_{i=1}^s b_i \left [\frac{\partial l}{\partial \eta}
(c_i h) \ddexp_{\ad_{\xi(c_i
h)}}\dot{\tilde{l}}_{\nu,s}(c_i)\right ]\, ,
\end{align*}
where $\nu=1,\ldots, s-1$.

\subsection{Discrete Euler--Poincar\'e Equations}
As shown above, we have constructed a higher-order reduced
discrete Lagrangian that depends on
\[ f_{kk+1}\equiv g_k g_{k+1}^{-1}.\]
We will now recall the derivation of the discrete Euler--Poincar\'e equations, introduced in \cite{MaPeSh1999}. The variations in $f_{kk+1}$ induced by variations in $g_k$, $g_{k+1}$ are computed as follows,
\begin{align*}
\delta f_{kk+1} &= -g_k^{-1} \delta g_k g_k{-1} g_{k+1} + g_k^{-1}\delta g_{k+1}\\
&= TR_{f_{kk+1}}(- g_k^{-1}\delta g_k + \Ad_{f_{kk+1}} g_{k+1}\delta g_{k+1})\, .
\end{align*}
Then, the variation in the discrete action sum is given by
\begin{align*}
\delta \mathbb{S}
&= \sum_{k=0}^{N-1} l'_d (f_{kk+1}) \delta f_{kk+1}\\
&= \sum_{k=0}^{N-1} l'_d (f_{kk+1}) TR_{f_{kk+1}}(- g_k^{-1}\delta g_k + \Ad_{f_{kk+1}} g_{k+1}\delta g_{k+1})\\
&= \sum_{k=1}^{N-1} \left [ l'_d(f_{k-1k})TR_{f_{k-1k}}\Ad_{f_{k-1k}} - l'_d(f_{kk+1})TR_{f_{kk+1}}\right ] \vartheta_k \, ,
\end{align*}
with variations of the form $\vartheta_k = g_k^{-1} \delta g_k$. In computing the variation of the discrete action sum, we have collected terms involving the same variations, and used the fact that $\vartheta_0=\vartheta_N=0$. This yields the {\bfi discrete Euler--Poincar\'e equation},
\begin{align*}
l'_d(f_{k-1k})TR_{f_{k-1k}}\Ad_{f_{k-1k}} -
l'_d(f_{kk+1})TR_{f_{kk+1}}&=0, & k&=1,\ldots, N-1.
\end{align*}
For ease of reference, we will recall the expressions from the
previous subsection that define the {\bfi higher-order reduced
discrete Lagrangian},
\[
l_d(f_{kk+1}) =h\sum_{i=1}^s b_i l\Big(\dexp_{\ad_{\xi(c_i
h)}}(\dot\xi(c_i h))\Big)\, ,
\]
where
\[ \xi(\xi^\nu;\tau h) = \sum_{\kappa=0}^s \xi^\kappa \tilde{l}_{\kappa,s}(\tau)\, ,\]
and
\begin{align*}
\xi^0 &= 0\, ,\\
\xi^s &= \log(f_{kk+1})\, ,
\end{align*}
and the remaining Lie algebra elements $\{\xi^\nu\}_{\nu=1}^{s-1}$, are defined implicitly by
\begin{align*}
0 &= h\sum_{i=1}^s b_i \left [\frac{\partial l}{\partial \eta}
(c_i h) \ddexp_{\ad_{\xi(c_i
h)}}\dot{\tilde{l}}_{\nu,s}(c_i)\right ],
\end{align*}
for $\nu=1,\ldots,s-1$, and where
\[
\ddexp_w = \sum_{n=0}^\infty \frac{w^n}{(n+2)!}\, .
\]
When the discrete Euler--Poincar\'e equation is used in
conjunction with the higher-order reduced discrete Lagrangian, we
obtain the {\bfi higher-order Euler--Poincar\'e equations}.

\section{Higher-Order Symplectic-Energy-Momentum Variational
Integrators}\label{gvi:sec:higher_order_sym_eng_mom}

In this section, we will generalize our new derivation of the
symplectic-energy-momentum preserving variational integrators
(see, \cite{KaMaOr1999}) to yield integrators with higher-order accuracy.

As before, we consider a piecewise interpolation, with both the
control points in the field variables, $q_i^\nu$, and the
endpoints of the interval, $h_i^0,h_i^1$, as degrees of freedom.
The continuity conditions for this function space are $q_i^s
=q_{i+1}^0$, and $h_i^1= h_{i+1}^0$. Then, we have that the
discrete action is given by
\begin{align*} \mathbb{S}_d =& \sum_{i=0}^{N-1} (h_i^1-h_i^0) \sum_{j=1}^s b_j L (q_i(c_j(h_i^1-h_i^0); q_i^\nu),\dot q_i(c_j(h_i^1-h_i^0); q_i^\nu))\\
&\qquad-\sum_{i=0}^{N-2} \lambda_i (q_i^s-q_{i+1}^0) -
\sum_{i=0}^{N-2} \omega_i(h_i^1-h_{i+1}^0),
\end{align*}
where
\begin{align*}
q_i(\tau (h_i^1-h_i^0);q_i^\nu)&=\sum_{\kappa=0}^s q_i^\kappa
\tilde{l}_{\kappa,s}(\tau),\\
\dot q_i(\tau
(h_i^1-h_i^0);q_i^\nu)&=\frac{1}{h_i^1-h_i^0}\sum_{\kappa=0}^s
q_i^\kappa \dot{\tilde{l}}_{\kappa,s}(\tau).
\end{align*}
To simplify the expressions, we define $h_i\equiv h_i^1-h_i^0$.
Then, the variational equations are given by
\begin{align*}
0 =& \sum_{j=1}^s b_j L(q_i(c_j h_i),\dot q_i(c_j h_i)) - h_i \sum_{j=1}^s b_j \frac{\partial L}{\partial \dot q}(c_j h_i) \frac{1}{h_i}\dot q_i(c_j h_i)-\omega_i, & \text{for } i=0,\ldots, N-2,\\
0 =& -\sum_{j=1}^s b_j L(q_i(c_j h_i),\dot q_i(c_j h_i)) + h_i \sum_{j=1}^s b_j\frac{\partial L}{\partial \dot q}(c_j h_i) \frac{1}{h_i}\dot q_i(c_j h_i)+\omega_{i-1}, &\text{for } i=1,\ldots, N-1,\\
0 =& h_i\sum_{j=1}^s b_j \left[\frac{\partial L}{\partial q}(c_j h_i) \tilde{l}_{s,s}(c_j)+\frac{1}{h_i}\frac{\partial L}{\partial \dot q}(c_j h_i)\dot{\tilde{l}}_{s,s}(c_j)\right] - \lambda_i, &\text{for }i=0,\ldots, N-2,\\
0 =&h_i\sum_{j=1}^s b_j \left[\frac{\partial L}{\partial q}(c_j h_i) \tilde{l}_{0,s}(c_j)+\frac{1}{h_i}\frac{\partial L}{\partial \dot q}(c_j h_i)\dot{\tilde{l}}_{0,s}(c_j)\right] + \lambda_{i-1}, &\text{for }i=1,\ldots, N-1,\\
0 =& h_i\sum_{j=1}^s b_j \left[\frac{\partial L}{\partial q}(c_j h_i) \tilde{l}_{\nu,s}(c_j)+\frac{1}{h_i}\frac{\partial L}{\partial \dot q}(c_j h_i)\dot{\tilde{l}}_{\nu,s}(c_j)\right],  &\stackrel{\displaystyle\text{for } i=0,\ldots, N-1,}{\hspace*{\fill}\nu=1,\ldots, s-1,}\\
q_i^s =& q_{i+1}^0, & \text{for } i=0,\ldots, N-2,\\
h_i^1=& h_{i+1}^0, & \text{for } i=0,\ldots, N-2.
\end{align*}
We can eliminate the Lagrange multipliers, to yield
\begin{align*}
0 =& \sum_{j=1}^s b_j L(q_i(c_j h_i),\dot q_i(c_j h_i))- \sum_{j=1}^s b_j\frac{\partial L}{\partial \dot q}(c_j h_i) \dot q_i(c_j h_i)\\
&\quad+\sum_{j=1}^s b_j L(q_{i-1}(c_j h_{i-1}),\dot q_{i-1}(c_j
h_{i-1}))\\
&\qquad- \sum_{j=1}^s b_j\frac{\partial L}{\partial \dot q}(c_j h_{i-1}) \dot q_{i-1}(c_j h_{i-1}), &\text{for }i=1,\ldots, N-1,\\
0 =& h_i\sum_{j=1}^s b_j \left[\frac{\partial L}{\partial q}(c_j h_i) \tilde{l}_{s,s}(c_j)+\frac{1}{h_i}\frac{\partial L}{\partial \dot q}(c_j h_i)\dot{\tilde{l}}_{s,s}(c_j)\right]\\
&\quad+ h_{i-1}\sum_{j=1}^s b_j \left[\frac{\partial L}{\partial q}(c_j h_{i-1}) \tilde{l}_{0,s}(c_j)+\frac{1}{h_{i-1}}\frac{\partial L}{\partial \dot q}(c_j h_{i-1})\dot{\tilde{l}}_{0,s}(c_j)\right], & \text{for }i=1,\ldots, N-1,\\
0 =& h_i\sum_{j=1}^s b_j \left[\frac{\partial L}{\partial q}(c_j h_i) \tilde{l}_{\nu,s}(c_j)+\frac{1}{h_i}\frac{\partial L}{\partial \dot q}(c_j h_i)\dot{\tilde{l}}_{\nu,s}(c_j)\right], &\stackrel{\displaystyle\text{for } i=0,\ldots, N-1,}{\hspace*{\fill}\nu=1,\ldots, s-1,}\\
q_i^s=&q_{i+1}^0,& \text{for } i=0,\ldots, N-2,\\
h_i^1=&h_{i+1}^0,& \text{for } i=0,\ldots, N-2.
\end{align*}
If we define the discrete Lagrangian as follows,
\[ L_d(q_i,q_{i+1},h_i) \equiv h_i \sum_{j=1}^s b_j L(q_i(c_j h_i),\dot q_i(c_j h_i)),\]
where
\begin{align*}
q_i(\tau h_i;q_i^\nu)&=\sum_{\kappa=0}^s q_i^\kappa \tilde{l}_{\kappa,s}(\tau),\\
\dot q_i(\tau h_i;q_i^\nu)&=\frac{1}{h_i}\sum_{\kappa=0}^s
q_i^\kappa \dot{\tilde{l}}_{\kappa,s}(\tau),
\end{align*}
and $q_i^0=q_i$, $q_i^s=q_1$, and the remaining terms were
defined implicitly by
\[
0 = h_i\sum_{j=1}^s b_j \left[\frac{\partial L}{\partial q}(c_j
h_i) \tilde{l}_{\nu,s}(c_j)+\frac{1}{h_i}\frac{\partial
L}{\partial \dot q}(c_j h_i)\dot{\tilde{l}}_{\nu,s}(c_j)\right],
\]
then the equations reduce to the following,
\begin{align*}
E_d(q_i,q_{i+1},h_i) &= -\frac{\partial}{\partial h_i} [ L_d(q_i,q_{i+1},h_i)],\\
E_d(q_i,q_{i+1},h_i) &= E_d(q_{i+1},q_{i+2},h_{i+1}),\\
0&=D_2 L_d(q_i, q_{i+1}, h_i) + D_1 L_d(q_{i+1},q_{i+2},h_{i+1}).
\end{align*}
which is a {\bfi higher-order symplectic-energy-momentum variational integrator}.

\paragraph{Solvability of the Energy Equation.} It should be noted that the discrete energy conservation equation is not necessarily solvable, in general, particularly near stationary points. This issue is discussed in \cite{KaMaOr1999, LeMaOrWe2004}, and can be addressed by reformulating the discrete energy conservation equation as an optimization problem that chooses the time step by minimizing the discrete energy error squared. Clearly, reformulating the discrete energy conservation equation yields the desired behavior whenever the discrete energy conservation equation can be solved, while allowing the computation to proceed when discrete energy conservation cannot be achieved, albeit with a slight energy error in that case. This does not degrade performance significantly, since instances in which discrete energy conservation cannot be achieved are rare.

\section{Spatio-Temporally Adaptive Variational Integrators}\index{variational integrator!spatio-temporally adaptive}
As is the case with all inner approximation techniques in numerical analysis, the quality of the numerical solution we obtain is dependent on the rate at which the sequence of finite-dimensional function spaces  approximates the actual solution as the number of degrees of freedom is increased.

For problems that exhibit shocks, nonlinear approximation spaces (see, for example, \cite{DeVore1998}), as opposed to linear approximation spaces, are clearly preferable. Adaptive techniques have been developed in the context of finite elements under the name of $r$-adaptivity and moving finite elements (see, for example, \cite{Ba1995}), and has been developed in a variational context for elasticity in \cite{ThOr2003}. The standard motivation in discrete mechanics to introduce function spaces that have degrees of freedom associated with the base space is to achieve energy or momentum conservation, as discussed in \S\ref{gvi:sec:higher_order_sym_eng_mom}, or \cite{KaMaOr1999,OlWeWu2004}. However, if the solution to be approximated exhibits shocks, nonlinear approximation techniques achieve better results for a given number of degrees of freedom.

In this section, we will sketch the use of regularizing transformations of the base space, to achieve a computational representation of sections of the configuration bundle that will yield more accurate numerical results.

\begin{figure}[htbp]
\begin{center}
\subfigure[Using equispaced nodes]{\includegraphics[scale=0.7]{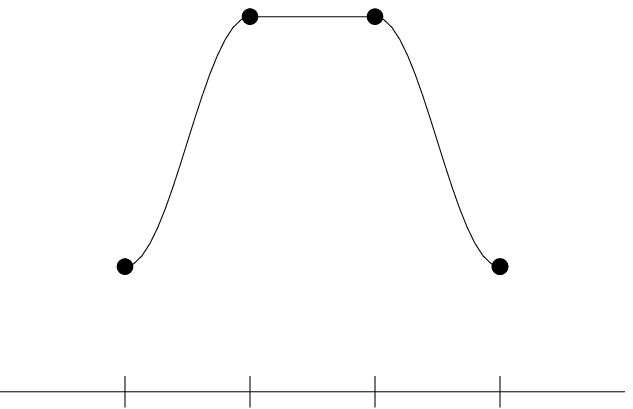}}
\subfigure[Using adaptive nodes]{\includegraphics[scale=0.7]{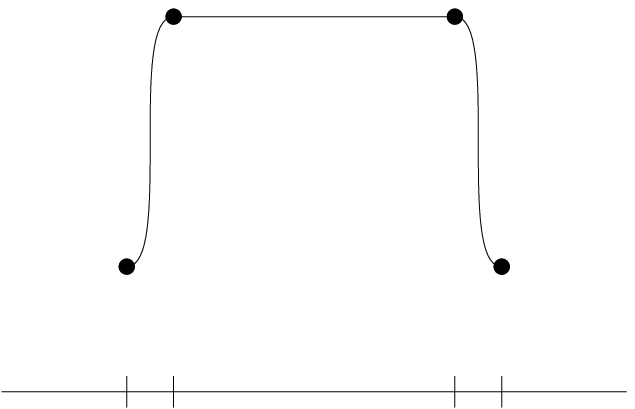}}
\end{center}
\caption{Linear and nonlinear approximation of a characteristic function.}\label{gvi:fig:linear_nonlinear_char_fn}
\end{figure}

Consider the situation when we are representing a characteristic function using piecewise spline interpolation. We show in Figure~\ref{gvi:fig:linear_nonlinear_char_fn}, the difference between linear and nonlinear approximation of the characteristic function.

When the derivatives of the solution vary substantially in a
spatially distributed manner, we obtain additional accuracy, for a
fixed representation cost, if we allow nodal points to cluster
near regions of high curvature. It is therefore desirable to
consider variational integrators based on function spaces that are
parameterized by both the position of the nodal points on the base
space, as well as the field values over the nodes.

\begin{figure}[htbp]
\begin{center}
\includegraphics{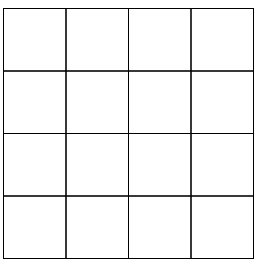}\raisebox{12mm}{$\quad\mapsto\quad$}
\includegraphics{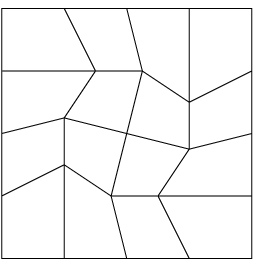}
\end{center}
\caption{Mapping of the base space from the computational to the physical domain.}\label{gvi:fig:computational_physical}
\end{figure}

This is represented by having a regular grid for the computational domain $\mathcal{R}$, which is then mapped to the physical base space $\mathcal{X}$, as shown in Figure~\ref{gvi:fig:computational_physical}.

\begin{figure}[htbp]
\[
\xymatrix{
\includegraphics[scale=0.7]{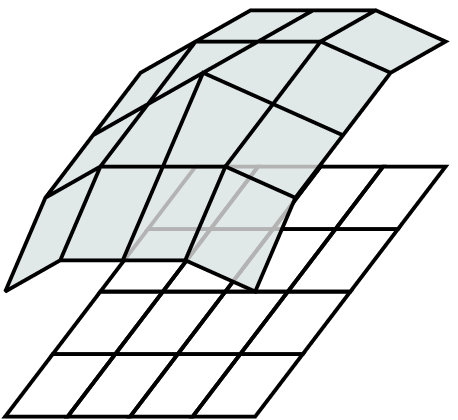} & \includegraphics[scale=0.7]{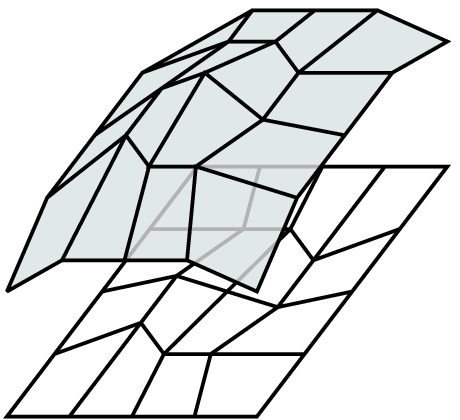}\\
\includegraphics{mesh1}\ar[r]_\varphi\ar[u]^{\tilde q} & \includegraphics{mesh2}\ar[u]_{q}
}\]
\caption{Factoring the discrete section.}\label{gvi:fig:section_factored}
\end{figure}

The sections of the configuration bundle factor as follows,
\[\xymatrix{
& Y\\
\mathcal{R}\ar[r]_\varphi \ar[ur]^{\tilde{q}}& \mathcal{X} \ar[u]_{q}
}
\]
The mapping $\varphi:\mathcal{R}\rightarrow \mathcal{X}$ results in a regularized computational representation $\tilde{q}:\mathcal{R}\rightarrow Y$ of the original section $q:\mathcal{X}\rightarrow Y$. The relationship between the discrete section of the configuration bundle and its computational representation is illustrated in Figure~\ref{gvi:fig:section_factored}.

The action integral is then given by
\[ \mathcal{S}(q)=\int_{\mathcal{X}}L(j^1 q)=\int_{\mathcal{R}}L(j^1 \tilde q) |\mathbf{D}\varphi|.
\]
Thus, even though $q:\mathcal{X}\rightarrow Y$ may exhibit shocks, the computational representation we work with, $\tilde{q}:\mathcal{R}\rightarrow Y$, is substantially more regular, and consequently, a numerical quadrature scheme in $\mathcal{R}$ applied to
\[ \int_{\mathcal{R}}L(j^1 \tilde q) |\mathbf{D}\varphi|\]
is significantly more accurate than the corresponding numerical quadrature scheme in $\mathcal{X}$ applied to
\[\int_{\mathcal{X}}L(j^1 q) .\]
As such, spatio-temporally adaptive variational integrators achieve increased accuracy by allowing the accurate representation of shock solutions using an adapted free knot representation, while using a smooth computational representation to compute the action integral.
%%%%%%%%%%%%%%%%%%%%%%%%%%%%%%%%%%%%%%%%%%%%%%
\section{Multiscale Variational Integrators}\label{gvi:sec:mvi}\index{variational integrator!multiscale}\index{multiscale!variational integrator|see{variational integrator, multiscale}}
In the work on multiscale finite elements (MsFEM), introduced and developed in \cite{HoWu1999, EfHoWu2000, ChHo2003}, shape functions that are solutions of the fast dynamics in the absence of slow forces are constructed to yield finite element schemes that achieve convergence rates that are independent of the ratio of fast to slow scales.

In constructing a multiscale variational integrator, we need to choose finite-dimensional function spaces that do a good job of approximating the fast dynamics of the problem, when the slow variables are frozen. In addition, we require an appropriate choice of numerical quadrature scheme to be able to evaluate the action, which involves integrating a highly-oscillatory Lagrangian. In this section, we will discuss how to go about making such choices of function spaces and quadrature methods.

We will start with a discussion of the multiscale finite element method, to illustrate the importance of a good choice of shape functions in computing solutions to problems with multiple scales. After that, we will walk through the construction of a multiscale variational integrator for the case of a planar pendulum with a stiff spring. Finally, we will discuss how we might proceed if we do not possess knowledge of which variables, or forces, are fast or slow.

\subsection{Multiscale Shape Functions}\index{multiscale!shape function}
We will illustrate the idea of constructing shape functions that are solutions of the fast dynamics by introducing a model multiscale second-order elliptic partial differential equation given by
\[ \nabla\cdot a(x/\epsilon)\nabla u^\epsilon(x)=f(x),\]
with homogeneous boundary conditions. In the one-dimensional case, we can solve for the solution analytically, and it has the form
\[ u^\epsilon(x)=\int_0^x \frac{F(y)}{a(y/\epsilon)}dy -\frac{\int_0^1 \frac{F(y)}{a(y/\epsilon)}dy}{\int_0^1\frac{dy}{a(y/\epsilon)}}\int_0^x \frac{dy}{a(y/\epsilon)},
\]
where $F(x)=\int_0^x f(y)dy$. If we have nodal points at $\{x_i\}_{i=0}^N$, then the appropriate multiscale shape functions to adopt in this example is to use shape functions that are solutions of the homogeneous problem at the element level. These shape functions $\varphi^\epsilon_i$ satisfy
\[\begin{cases}
\frac{\partial}{\partial x}\left(a(x/\epsilon)\frac{\partial}{\partial x}\varphi^\epsilon_i\right)=0, \quad \text{for } x_{i-1}<x<x_{i+1} ;\\
\varphi^\epsilon_i(x_{i-1}) =0;\quad
\varphi^\epsilon_i(x_{i+1}) =0;\quad
\varphi^\epsilon_i(x_i) =1.
\end{cases}
\]
And they have the explicit form given by
\[\varphi^\epsilon_i(x)=
\begin{cases}
\left[\int_{x_{i-1}}^{x_i} \frac{ds}{a(s/\epsilon)}\right]^{-1}\left[\int_{x_{i-1}}^x\frac{ds}{a(s/\epsilon)}\right], & x\in[x_{i-1},x_i]\, ;\\
\left[\int_{x_i}^{x_{i+1}}\! \frac{ds}{a(s/\epsilon)}\right]^{-1}\left[\int_x^{x_{i+1}}\!\frac{ds}{a(s/\epsilon)}\right], & x\in(x_i,x_{i+1}]\, ;\\
0, & \text{otherwise}\, .
\end{cases}\]
It can be shown that this will yield a numerical scheme that solves exactly for the solution at the nodal points. 

\begin{figure}[htbp]
\begin{center}
\subfigure[Exact solution]{\label{gvi:fig:ms1}\includegraphics[scale=0.7]{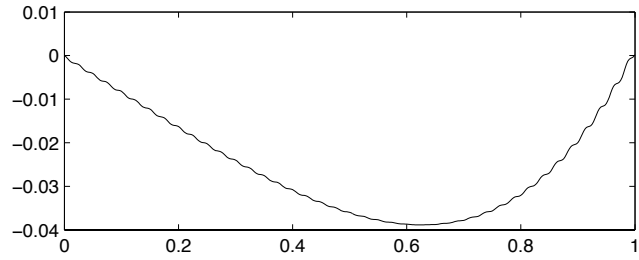}}\\
\subfigure[Exact solution (zoomed)]{\label{gvi:fig:ms2}\includegraphics[scale=0.7]{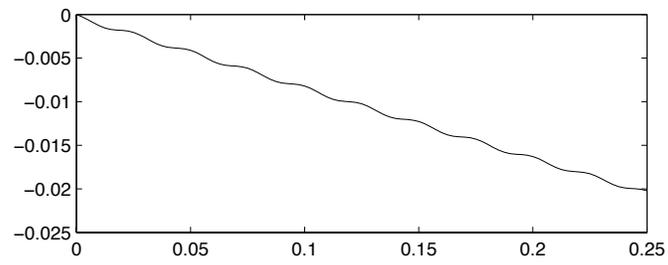}}\\
\subfigure[Multiscale shape function]{\label{gvi:fig:ms3}\includegraphics[scale=0.7]{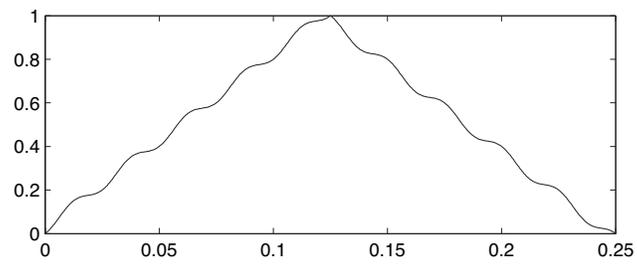}}
\caption{Comparison of the multiscale shape function and the exact solution for the elliptic problem.}
\label{gvi:fig:ms}
\end{center}
\end{figure}

As an example, we will compute the analytical solution for $a(x)=\frac{10}{1+0.95\sin(2\pi x)}$, $f(x)=x^2$, and $\epsilon=0.025$. This is illustrated in Figure~\ref{gvi:fig:ms1}. What is particularly interesting is to compare the zoomed plot of the exact solution and the multiscale shape function over the same interval, shown in Figures~\ref{gvi:fig:ms2} and \ref{gvi:fig:ms3}, respectively.
The multiscale finite element method is able to achieve excellent results because the multiscale shape functions are able to capture the fast dynamics well. In the next subsection, we will discuss how this insight is relevant in the construction of multiscale variational integrators.

\subsection{Multiscale Variational Integrator for the Planar Pendulum with a Stiff Spring}
\begin{figure}[htbp]
\begin{center}
\includegraphics[scale=0.15]{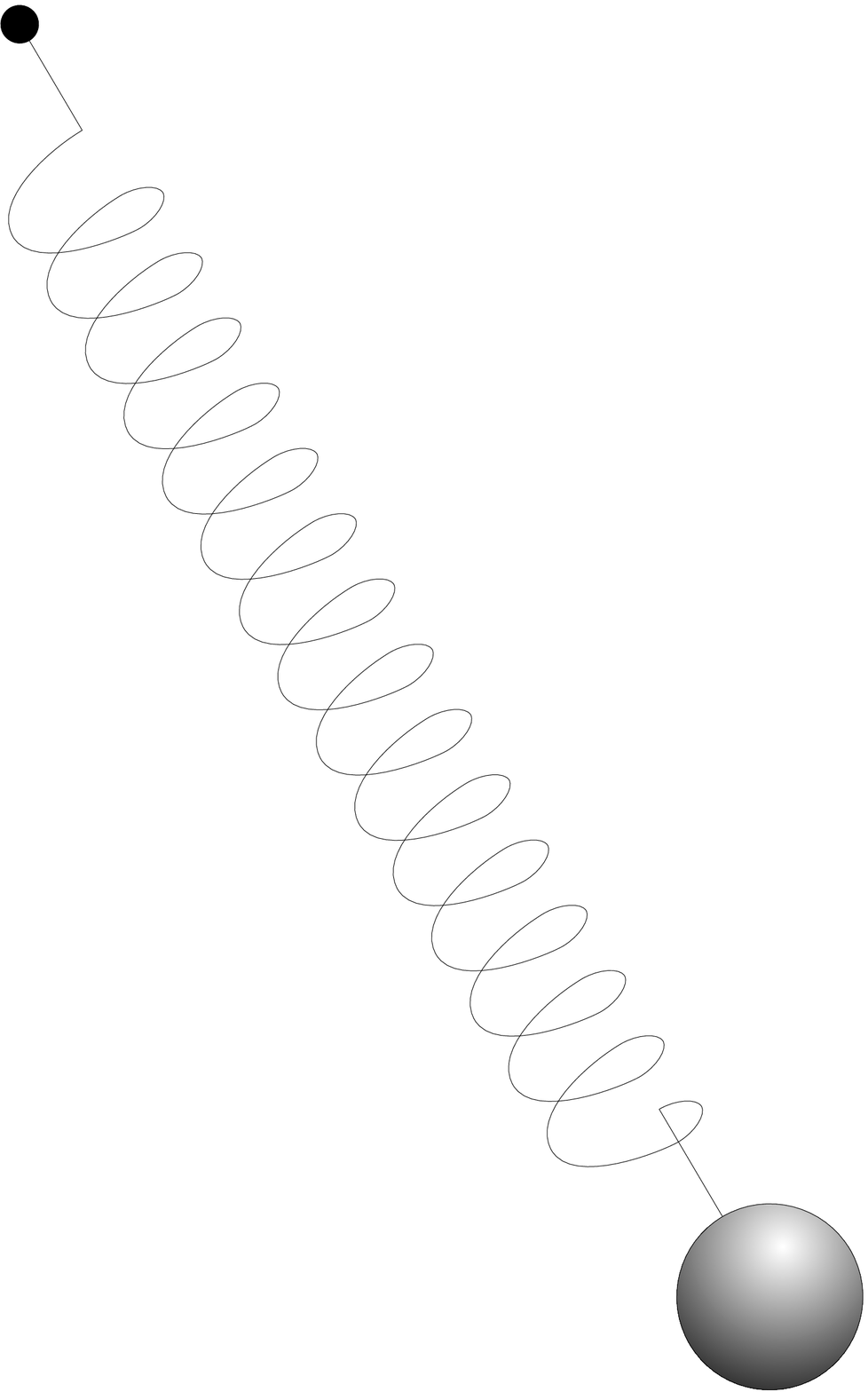}
\caption{Planar pendulum with a stiff spring}
\label{gvi:fig:stiff_pendulum}
\end{center}
\end{figure}

\index{variational integrator!multiscale!planar pendulum with a stiff spring}
As was shown previously, a shape function that captures the fast dynamics of a multiscale problem is able to achieve superior accuracy when used for computation. While this idea has primarily been used for problems with multiple spatial scales, it is natural to consider its application to a problem with multiple temporal scales, such as the problem of the planar pendulum with a stiff spring, as illustrated in Figure~\ref{gvi:fig:stiff_pendulum}.

We will use this example to illustrate the issues that arise in constructing a multiscale variational integrator. The variables are $q=(a,\theta)$, where $a$ is the spring extension, and $\theta$ is the angle from the vertical. The Lagrangian is given by
\[ L(a,\theta,\dot a,\dot\theta)=\frac{m}{2}(\dot x^2 +\dot y^2)-mgy-\frac{k}{2}a^2,\]
where
\begin{align*}
x &= (l+a)\sin\theta,\\
y &= -(l+a)\cos\theta,\\
\dot x &= \dot a \sin\theta +\dot\theta(l+a)\cos\theta,\\
\dot y &= -\dot a \cos\theta +\dot\theta(l+a)\sin\theta.
\end{align*}

The Hamilton's equations for the planar pendulum with a stiff
spring are
\begin{align*}
\dot a &= \frac{p_a}{m},\\
\dot \theta &= \frac{p_\theta}{m(l+a)^2},\\
\dot p_a &= -ka +g m \cos \theta + m(l+a)\dot\theta^2,\\
\dot p_\theta &= -gm (l+a)\sin\theta.
\end{align*}
The timescale arising from the mass-spring system is
$2\pi\sqrt{{m}/{k}}$. The timescale arising from the planar
pendulum system is $2\pi\sqrt{{l}/{g}}$. The ratio of timescales is given by $\epsilon=\sqrt{mg/kl}$.

\begin{figure}[htbp]
\begin{center}
\subfigure[$a(t)$]{\includegraphics[scale=0.7]{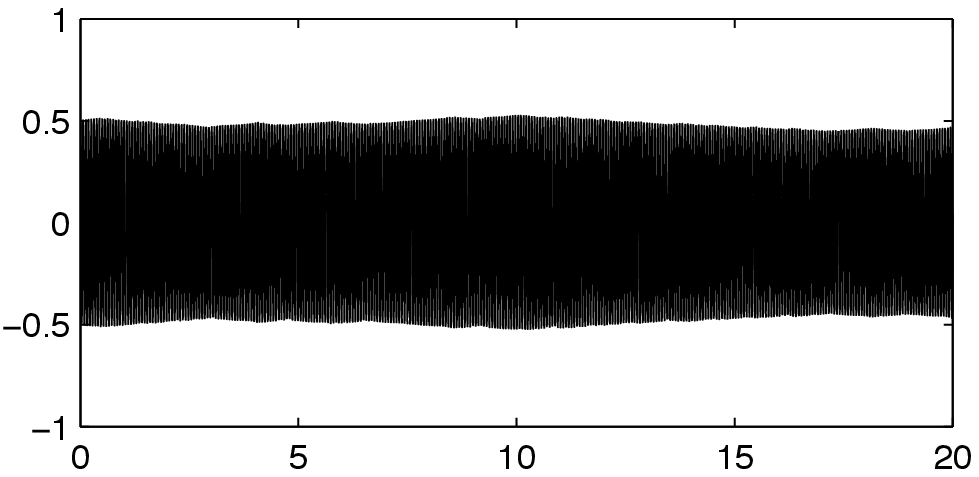}}\quad
\subfigure[$\theta(t)$]{\includegraphics[scale=0.7]{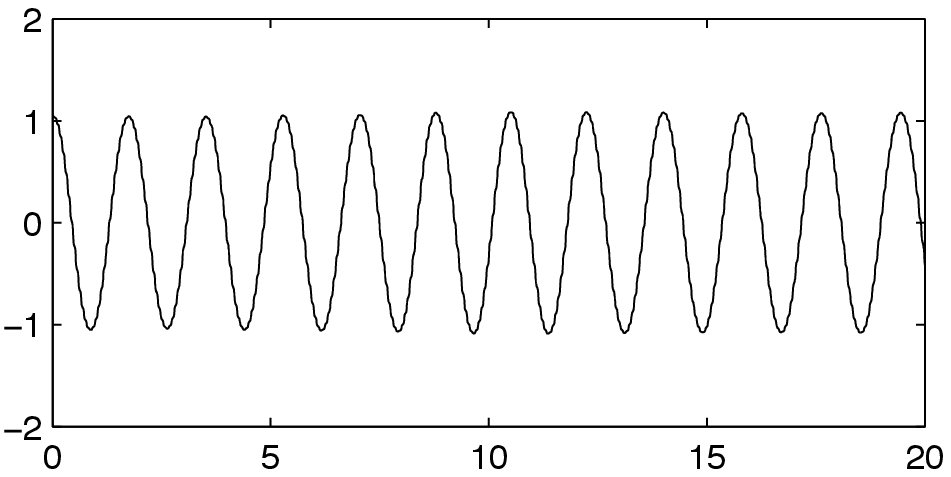}}\\
\subfigure[$a(t)$ (zoomed)]{\includegraphics[scale=0.7]{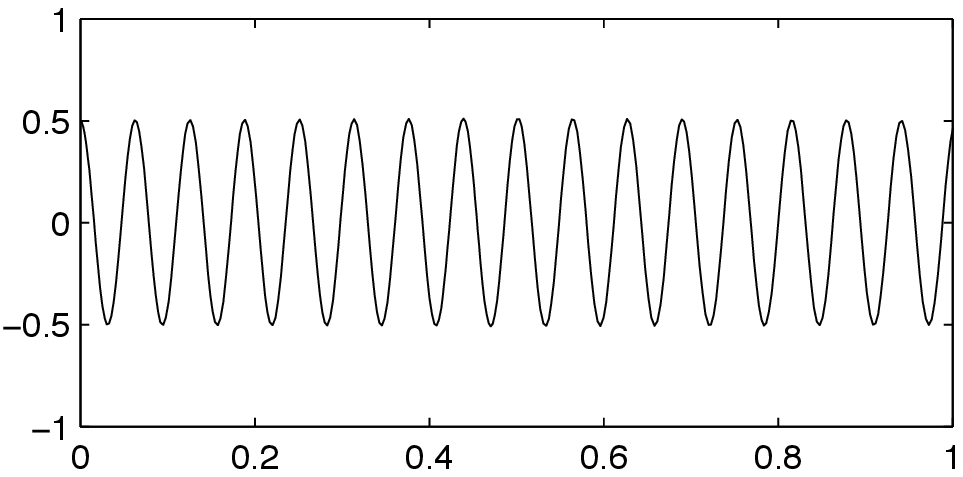}}\quad
\subfigure[$0.5\cos\left(\sqrt{k/m}\,t\right)$]{\includegraphics[scale=0.7]{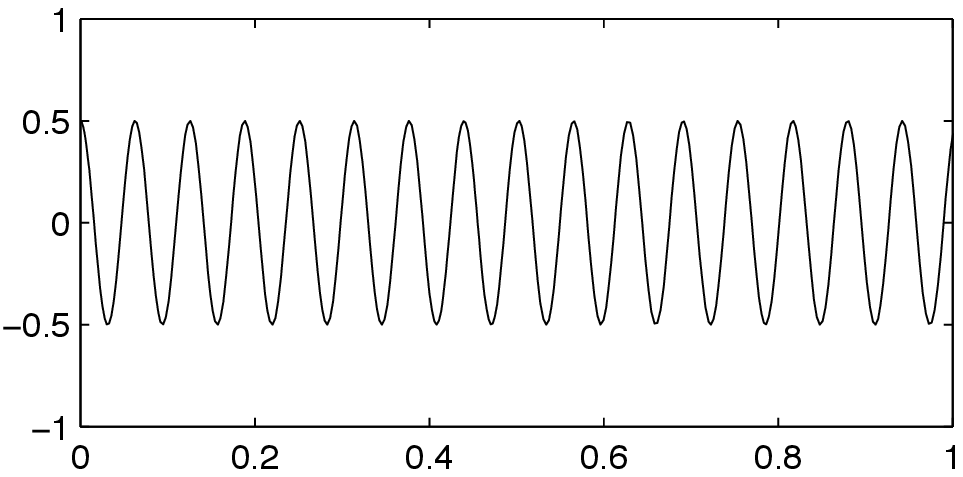}}\\
\caption{Comparison of the multiscale shape function and the exact solution for the planar pendulum with a stiff spring.}
\label{gvi:fig:mvi_spring}
\end{center}
\end{figure}

\paragraph{Multiscale Shape Function.}\index{multiscale!shape function}
In this problem, the fast scale is associated with the stiff spring, and if we set the slow variable $\theta=0$, we obtain the equation
\[\ddot a = \frac{\dot p_a}{m} = -\frac{k}{m}a,\]
which has solutions of the form
\[ a(t) = a_0\sin\left(\sqrt{k/m}\,t\right)+a_1\cos\left(\sqrt{k/m}\,t\right).\]
We will now consider a well-resolved simulation of this system using the \texttt{ode15s} stiff solver from \textsc{Matlab}, with parameters $m=1$, $g=9.81$, $k=10000$, $l=1$, giving a scale separation of $\epsilon=0.0313$. The simulation results are show in Figure~\ref{gvi:fig:mvi_spring}.

Clearly, if we wish to choose time steps that do not resolve the fast oscillations in $a$, but do resolve the slow oscillations in $\theta$, it would be desirable to include $\sin\left(\sqrt{k/m}\,t\right)$, and $\cos\left(\sqrt{k/m}\,t\right)$ in the finite-dimensional function space used to interpolate $a(t)$.

\paragraph{Evaluating the Discrete Lagrangian.}\index{highly-oscillatory integral}
Since we have chosen time steps that do not resolve the fast oscillations, it follows that over the interval $[0,h]$, the Lagrangian will oscillate rapidly as well. In computing the discrete Lagrangian, it is therefore necessary to ensure that this highly-oscillatory integral is well-approximated.

It is conventional wisdom, in numerical analysis, that the numerical quadrature of highly-oscillatory integrals is a challenging problem requiring the use of many function evaluations. Recently, there has been a series of papers, \cite{Is2003a, Is2003b, Is2004, IsNo2004}, that provides an analysis of Filon-type quadrature schemes that provide an efficient and accurate method of evaluating such integrals. The method is applicable to weighted integrals as well, but we will summarize the results from \S 3 of \cite{Is2003a} restricted to unweighted integrals, and refer the reader to the original reference for an in-depth discussion and analysis.\index{Filon quadrature}

The Filon-type method aims to evaluate an integral of the form
\[ I_h[f] = \int_0^h f(x) e^{i\omega x} dx = h \int_0^h f(hx) e^{i\omega x} dx.\]
Given a set of distinct quadrature points, $c_1<c_2<\cdots<c_\nu$ in $[0,1]$, the Filon-type quadrature method is given by
\[ Q^F_h[f]=h\sum_{i=1}^\nu b_i (ih\omega) f (c_i h),\]
where
\[ b_i(ih\omega)=\int_0^1 l_i(x) e^{ih\omega x} dx,\]
and $l_i$ are the Legendre polynomials. Here, we draw attention to the fact that the quadrature weights are dependent on $h\omega$.

If the quadrature points correspond to Gauss-Christoffel quadrature of order $p$, then the error for the Filon-type method is given as follows.
\begin{align*}
\mathcal{O}(h^{p+1}),\quad&\text{if}\quad h\omega\ll 1\\
\mathcal{O}(h^\nu),\quad&\text{if}\quad h\omega=\mathcal{O}(1)\\
\mathcal{O}(h^{\nu+1}/(h\omega)),\quad&\text{if}\quad h\omega\gg
1\\
\mathcal{O}(h^{\nu+1}/(h\omega)^2),\quad&\text{if}\quad h\omega\gg
1, c_1=0, c_\nu=1.
\end{align*}
Clearly, for highly-oscillatory functions, the last case, which corresponds to the Lobatto quadrature points, is most desirable. Thus, it is appropriate to use the Filon-Lobatto method to evaluate the discrete Lagrangian in our case.

\paragraph{Discrete Variational Equations.}
As we discussed previously, instead of looking for stationary solutions to the discrete Hamilton's principle in polynomial spaces, we will consider solutions that are piecewise of the form
\[q(t;\{p_j\},\omega,a_0,a_1) = \left(\sum\nolimits_{j=0}^n p_j t^j\right)
(1+a_0\sin(\omega t)+a_1\cos(\omega t))\, .\]

This function space approximates the highly-oscillatory nature of the solution well, in contrast to a polynomial function space, thereby avoiding approximation-theoretic errors. The degrees of freedom in this function space are $\{p_j\}$, $\omega$, $a_0$, and $a_1$. Since there are no distinguished degrees of freedom that are responsible for the endpoint values of the curve, we need to impose continuity at the nodes using a Lagrange multiplier. The augmented discrete action is given by
\begin{align*}
\mathbb{S}_d =& \sum_{i=0}^{N-1} \int_0^h L(j^1 q_i(t;\{p_j^i\},\omega^i,a_0^i,a_1^i) ) dt\\*
&\qquad+\sum_{i=0}^{N-2} \lambda_i (q_i(h;\{p_j^i\},\omega^i,a_0^i,a_1^i)-q_{i+1}(0;\{p_j^{i+1}\},\omega^{i+1},a_0^{i+1},a_1^{i+1}))\, ,
\end{align*}
where each of the integrals are evaluated using the Filon-Lobatto method. Taking variations with respect to the degrees of freedom yields an update map,
\[\left(\{p_j^i\},\omega^i,a_0^i,a_1^i\right)\mapsto \left(\{p_j^{i+1}\},\omega^{i+1},a_0^{i+1},a_1^{i+1}\right),\]
which gives the {\bfi multiscale variational integrator} for the planar pendulum with a stiff spring.

\subsection{Computational Aspects}
Multiscale variational integrators have the advantage of directly accounting for the contribution of the fast dynamics, thereby allowing the scheme to use significantly larger time-steps, while maintaining accuracy and stability. It is possible to take advantage of knowledge about which of the variables, or forces, are fast or slow, by using a low degree polynomial and oscillatory functions for the fast variables, and a higher-order polynomial for the slow variables. In the absence of such information, it is appropriate to use a function space with both polynomials and oscillatory functions, and apply it to all the variables.

Recall that the Filon-type method has quadrature coefficients that depend on the frequency. As such, the initial fast frequency has to be estimated numerically using a fully resolved computation for a short period of time. Since both the function space and the
quadrature weights depend on the fast frequency $\omega$, the resulting scheme is implicit and fairly nonlinear, and as such, it may be expensive for large systems.

%%%%%%%%%%%%%%%%%%%%%%%%%%%%%%%%%%%%%%%%%%%%%%
\section{Pseudospectral Variational Integrators}\index{variational integrator!pseudospectral}
The use of spectral expansions of the solution in space are particularly
appropriate for highly accurate simulations of the evolution of
smooth solutions, such as those arising from quantum mechanics. We
will introduce pseudospectral variational integrators, and consider the
Schr\"odinger equation as an example.

In particular we will adopt the tensor product of a spectral expansion in space, and a polynomial expansion in time. For example, we could have an interpolatory function of the form
\[\psi(x,(\tau+l)\Delta t) =
\frac{1}{2\pi}\sideset{}{'}{\sum}_{k=-N/2}^{N/2} e^{ikx}
\left((1-\tau) \hat v_k^l +\tau\hat v_k^{l+1}\right),\]
which is the tensor product of a discrete Fourier expansion in space, and linear interpolation in time. Here, the $\sideset{}{'}{\textstyle\sum}$ notation denotes a weighted sum where the terms with indices $\pm N/2$ are weighted by $1/2$, and the other terms are weighted by $1$. See page 19 of \cite{Trefethen2000} for a discussion of why this is necessary to fix an issue with derivatives of the interpolant.

The degrees of freedom are given by $\hat v^l_k$, which are the discrete Fourier coefficients. We will later see how such an interpolation can be applied to the Schr\"odinger equation. The action integral can be exactly evaluated for this class of shape functions, as we will see below.

It is straightforward to generalize the pseudospectral approach we present in this section to a spectral variational integrator, with discrete Fourier expansions in space for periodic domains, or Chebyshev expansions in space for non-periodic domains, and Chebyshev expansions in time. This will however result in all the degrees of freedom on the space-time mesh being coupled, and is therefore substantially more expensive computationally than the pseudospectral method. The payoff for adopting the spectral approach is spectral accuracy, which is accuracy beyond all orders.

\subsection{Variational Derivation of the Schr\"{o}dinger
Equation}\index{variational problems!Schr\"odinger equation|see{Schr\"{o}dinger equation}}\index{Schr\"{o}dinger equation} Let $\mathcal{H}$ be a complex Hilbert space, for
example, the space of complex-valued functions $\psi$ on
$\mathbb{R}^3$ with the Hermitian inner product,
\[ \langle\psi_1,\psi_2\rangle=\int \psi_1(x)\overline{\psi}_2(x)d^3 x\, ,\]
where the overbar denotes complex conjugation. We will present a Lagrangian derivation of the Schr\"odinger equation, following worked example 9.1 on pages 568--569 of~\cite{JoSa1998}.

Consider the Lagrangian density $\mathcal{L}$ given by\index{Schr\"{o}dinger equation!action functional}
\[ \mathcal{L}(j^1 \psi)=\frac{i\hbar}{2}\{\dot\psi\overline{\psi}-\psi\dot{\overline\psi}\}-
\hat H \psi\overline{\psi},\] where $\hat H:\mathcal{H}\rightarrow
\mathcal{H}$ is given by
\[ \hat H\psi=-\frac{\hbar^2}{2m}\nabla^2\psi+V\psi,\]
which yields
\[ \mathcal{L}(j^1 \psi)=\frac{i\hbar}{2}\{\dot\psi\overline{\psi}-\psi\dot{\overline\psi}\}-\frac{\hbar^2}{2m}\nabla\psi\cdot\nabla\overline{\psi}-
V \psi\overline{\psi}.\]
We take $\psi,\overline{\psi}$ as independent variables, and compute,
\begin{align*}\delta\int L dt
&=\int \left[\left( \frac{\partial\mathcal{L}}{\partial
\overline{\psi}}\delta\overline{\psi} +
\frac{\partial\mathcal{L}}{\partial
\dot{\overline{\psi}}}\delta\dot{\overline{\psi}} +\frac{\partial
\mathcal{L}}{\partial \nabla
\overline{\psi}}\delta\nabla\overline{\psi}\right)+\Biggl(\frac{\partial\mathcal{L}}{\partial
\psi}\delta\psi + \frac{\partial\mathcal{L}}{\partial
\dot\psi}\delta\dot\psi +\frac{\partial \mathcal{L}}{\partial
\nabla \psi}\delta\nabla\psi\Biggr)\right]d^3 x
dt\\
&=\int \left[\left( \frac{\partial\mathcal{L}}{\partial
\overline{\psi}} - \frac{\partial}{\partial
t}\frac{\partial\mathcal{L}}{\partial \dot{\overline{\psi}}}
-\nabla\cdot\frac{\partial\mathcal{L}}{\partial
\nabla\overline{\psi}}\right)\delta\overline{\psi}+\Biggl(\frac{\partial\mathcal{L}}{\partial
\psi} - \frac{\partial}{\partial
t}\frac{\partial\mathcal{L}}{\partial \dot\psi}
-\nabla\cdot\frac{\partial\mathcal{L}}{\partial
\nabla\psi}\Biggr)\delta\psi \right]d^3 x dt\\&=\int \left[
\left(\frac{i\hbar}{2}\dot\psi-V\psi+\frac{i\hbar}{2}\dot\psi+\frac{\hbar^2}{2m}\nabla^2\psi\right)\delta\overline{\psi}+
\left(\frac{i\hbar}{2}\dot{\overline{\psi}}-V\overline{\psi}+\frac{i\hbar}{2}\dot{\overline{\psi}}+\frac{\hbar^2}{2m}\nabla^2\overline{\psi}\right)\delta\psi\right]d^3x
dt,\end{align*} where we integrated by parts, and neglected
boundary terms as the variations vanish at the boundary of the
space-time region. Since the variations are arbitrary, we obtain the nonrelativistic
(linear) Schr\"odinger equation as a result,
\[i\hbar \dot\psi = \left\{-\frac{\hbar^2}{2m}\nabla^2+V\right\}\psi\, .\]

We note that the Lagrangian density is invariant under the
internal phase shift given by
\[\psi\mapsto e^{i\epsilon}\psi,\qquad \overline{\psi}\mapsto
e^{-i\epsilon}\overline{\psi}.\]

The space part of the multi-momentum map is given by
\[j^k = \frac{\partial \mathcal{L}}{\partial
(\partial_k\psi)}i\psi+\frac{\partial \mathcal{L}}{\partial
(\partial_k\overline{\psi})}(-i\overline{\psi})=\frac{i\hbar^2}{2m}(\psi\partial_k\overline\psi-\overline\psi\partial_k\psi),\]
and the time part is given by
\[j^0=\frac{\partial\mathcal{L}}{\partial\dot\psi}i\psi-\frac{\partial\mathcal{L}}{\partial
\dot{\overline\psi}}(-i\overline\psi)
=-\frac{\hbar}{2}(\dot{\overline\psi}\psi-\dot\psi\overline\psi).\]

The norm of the wavefunction is automatically preserved by
variational integrators, since the norm is a quadratic invariant.

\subsection{Pseudospectral Variational Integrator for the Schr\"odinger Equation}\index{variational integrator!pseudospectral!Schr\"{o}dinger equation}
Consider a periodic domain $[0,2\pi]$, discretized with a discrete
Fourier series expansion in space, and a linear interpolation in time. Let $N$ be an even integer, then, our computation is done on the following mesh,
\[\begin{xy}0;<10mm,0cm>:<0cm,10mm>::
(0,0.3)*\txt{$0$}, (5,0.3)*\txt{$\pi$}, (10,0.3)*\txt{$2\pi$},
(1,-0.3)*\txt{$x_1$}, (2,-0.3)*\txt{$x_2$},
(5,-0.3)*\txt{$x_{N/2}$}, (9,-0.3)*\txt{$x_{N-1}$},
(10,-0.3)*\txt{$x_N$}, (1,0)*@{*}, (2,0)*@{*}, (3,0)*@{*},
(4,0)*@{*}, (5,0)*@{*}, (6,0)*@{*}, (7,0)*@{*}, (8,0)*@{*},
(9,0)*@{*}, (10,0)*@{*}, \ar@{-} (0,0);(10,0), \ar@{-}
(0,-0.1);(0,0.1)
\end{xy}
\]
This implies that the grid spacing is given by
\[ h=\frac{2\pi}{N}\,.\]
The interpolation is given by
\begin{align*}
\psi(x,(\tau+l)\Delta t) &=
\frac{1}{2\pi}\sideset{}{'}{\sum}_{k=-N/2}^{N/2} e^{ikx}
\left((1-\tau) \hat v_k^l +\tau\hat v_k^{l+1}\right),\\
\dot\psi(x,(\tau+l)\Delta t) &= \frac{1}{2\pi\Delta
t}\sideset{}{'}{\sum}_{k=-N/2}^{N/2} e^{ikx} \left(\hat v_k^{l+1}
- \hat v_k^l\right),\\
\bar\psi(x,(\tau+l)\Delta t) &=
\frac{1}{2\pi}\sideset{}{'}{\sum}_{k=-N/2}^{N/2} e^{-ikx}
\left((1-\tau) \vbar_k^l +\tau\vbar_k^{l+1}\right),\\
\dot{\bar\psi}(x,(\tau+l)\Delta t) &= \frac{1}{2\pi\Delta
t}\sideset{}{'}{\sum}_{k=-N/2}^{N/2} e^{-ikx} \left(\vbar_k^{l+1}
- \vbar_k^l\right),
\end{align*}
and the discrete Fourier transformation is given by
\[ \hat v_j = \frac{1}{2\pi}h \sum_{j=1}^N e^{-ikx_j} v_j\, ,
\]
for $k=-N/2+1,\ldots,N/2$, and $\hat v_{-N/2}\equiv \hat v_{N/2}$.
Recall that
\[ \mathcal{L}(j^1\psi)=\frac{i\hbar}{2}\{\dot\psi\overline{\psi}-\psi\dot{\overline\psi}\}-
\hat H \psi\overline{\psi},\] where $\hat H:\mathcal{H}\rightarrow
\mathcal{H}$ is given by
\[ \hat H\psi=-\frac{\hbar^2}{2m}\nabla^2\psi+V\psi.\]
Furthermore, the potential $V$ is expressed using a discrete Fourier expansion,
\[ V(x)=\frac{1}{2\pi}\sideset{}{'}{\sum}_{k=-N/2}^{N/2} e^{ikx}\hat V_k\, .\]
In addition, we will need to introduce a normalization condition, so as to eliminate trivial solutions of the partial differential equation. The normalization condition is
\begin{align*}
1=\langle\psi_l,\psi_l\rangle
=
\int_0^{2\pi}
\left(\frac{1}{2\pi}\sideset{}{'}{\sum}_{k=-N/2}^{N/2} e^{ikx}\hat v_k^l \right)
\left(\frac{1}{2\pi}\sideset{}{'}{\sum}_{k=-N/2}^{N/2} e^{-ikx}\vbar_k^l \right)dx
= \frac{1}{2\pi}\,\,\sideset{}{''}{\sum}_{k=-N/2}^{N/2}\hat v_k^l\vbar_k^l,
\end{align*}
which is enforced using a Lagrange multiplier.
\paragraph{Discrete Action for the Schr\"odinger Equation.}\index{Schr\"{o}dinger equation!discrete action}
The discrete action in the space-time region
$[0,2\pi]\times[l\Delta t,(l+1)\Delta t]$ is given by
\begin{align*}
\mathbb{S}_d
&= \int_{l\Delta t}^{(l+1)\Delta t}\int_0^{2\pi} \mathcal{L}(j^1\psi) dx dt+\lambda_l(1-\langle\psi_l,\psi_l\rangle)\\
&= \int_{l\Delta t}^{(l+1)\Delta t}\int_0^{2\pi} \left[\frac{i\hbar}{2}\{\dot\psi\bar\psi-\psi\dot{\bar\psi}\}+\frac{\hbar^2}{2m}\nabla^2\psi\bar\psi-V\psi\bar\psi\right]dxdt + \lambda_l\left(1-\frac{1}{2\pi}\,\,\sideset{}{''}{\sum}_{k=-N/2}^{N/2}\hat v_k^l\vbar_k^l\right)\\
&= \int_0^1 \int_0^{2\pi} \frac{i\hbar}{2} \left[
\left(\frac{1}{2\pi\Delta t}\sideset{}{'}{\sum}_{k=-N/2}^{N/2}
e^{ikx} \left(\hat v_k^{l+1} - \hat v_k^l\right) \right)
\left(\frac{1}{2\pi}\sideset{}{'}{\sum}_{k=-N/2}^{N/2} e^{-ikx}
\left((1-\tau) \vbar_k^l +\tau\vbar_k^{l+1}\right)\right)\right.\\*
&\qquad\qquad\left.-\left(
\frac{1}{2\pi}\sideset{}{'}{\sum}_{k=-N/2}^{N/2} e^{ikx}
\left((1-\tau) \hat v_k^l +\tau\hat v_k^{l+1}\right)\right)
\left(\frac{1}{2\pi\Delta t}\sideset{}{'}{\sum}_{k=-N/2}^{N/2}
e^{-ikx} \left(\vbar_k^{l+1} - \vbar_k^l\right) \right)\right]\Delta t\,dx d\tau\\
&\qquad+\int_0^1 \int_0^{2\pi}
\frac{\hbar^2}{2m}\left(\frac{1}{2\pi}\sideset{}{'}{\sum}_{k=-N/2}^{N/2}
(-k^2)e^{ikx} \left((1-\tau) \hat v_k^l +\tau\hat
v_k^{l+1}\right)\right)\\*
&\qquad\qquad\cdot
\left(\frac{1}{2\pi}\sideset{}{'}{\sum}_{k=-N/2}^{N/2} e^{-ikx}
\left((1-\tau) \vbar_k^l +\tau\vbar_k^{l+1}\right)\right)\Delta t\, dxd\tau\\
&\qquad-\int_0^1\int_0^{2\pi}\left[\left(\frac{1}{2\pi}\sideset{}{'}{\sum}_{k=-N/2}^{N/2}e^{ikx}\hat V_k\right)
\left( \frac{1}{2\pi}\sideset{}{'}{\sum}_{m=-N/2}^{N/2} e^{imx}
\left((1-\tau) \hat v_m^l +\tau\hat v_m^{l+1}\right)\right)\right.\\*
&\qquad\qquad\cdot
\left.\left(\frac{1}{2\pi}\sideset{}{'}{\sum}_{n=-N/2}^{N/2} e^{-inx}
\left((1-\tau) \vbar_n^l +\tau\vbar_n^{l+1}\right)\right)\right]\Delta t\, dx dt\\
&\qquad+ \lambda_l\left(1-\frac{1}{2\pi}\,\,\sideset{}{''}{\sum}_{k=-N/2}^{N/2}\hat v_k^l\vbar_k^l\right)\\
&= \int_0^1 \frac{i\hbar}{2} \left[\frac{1}{2\pi}\sideset{}{''}{\sum}_{k=-N/2}^{N/2} \Bigl((\hat v^{l+1}_k-\hat v^l_k)((1-\tau) \vbar^l_k+\tau \vbar^{l+1}_k)-((1-\tau) \hat v_k^l +\tau\hat v_k^{l+1}) (\vbar^{l+1}_k-\vbar^l_k)\Bigr) \right] d\tau\\*
&\qquad-\int_0^1\left[\frac{\hbar^2}{2\pi} \frac{k^2}{2\pi}
\sideset{}{''}{\sum}_{k=-N/2}^{N/2} ((1-\tau)\hat v^l_k + \tau
\hat v^{l+1}_k)((1-\tau)\vbar^l_k + \tau
\vbar^{l+1}_k)\right] \Delta t\, d\tau\\*
&\qquad-\int_0^1\left(\frac{1}{2\pi}\right)^2\left[
\sideset{}{'}{\sum}_{n=-N/2}^{-1}\,\sideset{}{'}{\sum}_{m=-N/2}^{N/2+n}
\Bigl(\hat V_{n-m} ((1-\tau)\hat v^l_m +\tau\hat v^{l+1}_m) ((1-\tau)\vbar^l_n +\tau\vbar^{l+1}_n)\Bigr)\right.\\*
&\qquad\qquad+\left.
\sideset{}{'}{\sum}_{n=0}^{N/2}\,\sideset{}{'}{\sum}_{m=n-N/2}^{N/2} \Bigl(\hat V_{n-m} ((1-\tau)\hat v^l_m +\tau\hat v^{l+1}_m) ((1-\tau)\vbar^l_n +\tau\vbar^{l+1}_n)\Bigr)
\right]\Delta t\, d\tau\\*
&\qquad+ \lambda_l\left(1-\frac{1}{2\pi}\,\,\sideset{}{''}{\sum}_{k=-N/2}^{N/2}\hat v_k^l\vbar_k^l\right)\\
&= \frac{i\hbar}{4\pi}\,\,\sideset{}{''}{\sum}_{k=-N/2}^{N/2}\Bigl[\hat v^{l+1}_k \vbar^l_k - \hat v^l_k \vbar^{l+1}_k\Bigr]
-\frac{\hbar^2 k^2\Delta t}{24\pi^2} \sideset{}{''}{\sum}_{k=-N/2}^{N/2} \Bigl[\hat v^l_k (2\vbar^l_k+\vbar^{l+1}_k)+\hat v^{l+1}_k(\vbar^l_k+2\vbar^{l+1}_k)\Bigr]\\*
&\qquad-\frac{\Delta t}{24\pi^2}\left(\sideset{}{'}{\sum}_{n=-N/2}^{-1}\,\sideset{}{'}{\sum}_{m=-N/2}^{N/2+n} \,\,+\,\,\sideset{}{'}{\sum}_{n=0}^{N/2}\,\sideset{}{'}{\sum}_{m=n-N/2}^{N/2}\right)\hat V_{n-m}\Bigl[\hat v^l_m (2\vbar^l_n+\vbar^{l+1}_n)+\hat v^{l+1}_m(\vbar^l_n+2\vbar^{l+1}_n)\Bigr]\\*
&\qquad+ \lambda_l\left(1-\frac{1}{2\pi}\,\,\sideset{}{''}{\sum}_{k=-N/2}^{N/2}\hat v_k^l\vbar_k^l\right)\, ,
\end{align*}
where we used the fact that
\[
\int_0^{2\pi} e^{ikx}dx =2\pi\delta^i_0\, ,
\]
for $k\in\mathbb{Z}$, and we define $\sideset{}{'}{\textstyle\sum}$ as a weighted sum where the terms with  indices $\pm N/2$ are weighted by $1/2$, and $\sideset{}{''}{\textstyle\sum}$ as a weighted sum where the terms with indices $\pm N/2$ are weighted by $1/4$. We should note that using the same approach, it would be possible to exactly evaluate the action integral for the class of tensor product shape functions with a discrete Fourier expansion in space, and a polynomial expansion in time. In particular, a similar approach is valid in exactly evaluating the action integral when we use shape functions that are spectral in both space and time.

\paragraph{Discrete Euler--Lagrange Equations.}\index{Schr\"{o}dinger equation!Euler--Lagrange} We are now in a position to compute the discrete Euler--Lagrange equations associated with the Schr\"odinger equation when using a tensor product of a discrete Fourier expansion in space, and a linear interpolation in time.

The discrete variational equations are given by
\begin{align*}
0 &=\frac{i\hbar}{4\pi}\left[\vbar^{l-1}_j-\vbar^{l+1}_j\right]-\frac{\hbar^2 k^2 \Delta t}{24\pi^2}\left[\vbar^{l-1}_j+4\vbar^l_j+\vbar^{l+1}_j\right]\\*
&\qquad-\frac{\Delta t}{24\pi^2}\sideset{}{'}{\sum}_{n=-N/2}^{N/2+j}\hat V_{n-j}\left[\vbar^{l-1}_n+4\vbar^l_n+\vbar^{l+1}_n\right]-\frac{\lambda_l}{2\pi}\vbar^l_j\, , & \text{for } j=-N/2+1,\ldots,-1,\\
0 &=\frac{i\hbar}{4\pi}\left[\hat v^{l+1}_j-\hat v^{l-1}_j\right]-\frac{\hbar^2 k^2 \Delta t}{24\pi^2}\left[\hat v^{l-1}_j+4\hat v^l_j+\hat v^{l+1}_j\right]\\*
&\qquad-\frac{\Delta t}{24\pi^2}\sideset{}{'}{\sum}_{n=-N/2}^{N/2+j}\hat V_{j-n}\left[\hat v^{l-1}_n+4\hat v^l_n+\hat v^{l+1}_n\right]-\frac{\lambda_l}{2\pi}\hat v^l_j\, , & \text{for } j=-N/2+1,\ldots,-1,\\
0 &=\frac{i\hbar}{4\pi}\left[\vbar^{l-1}_j-\vbar^{l+1}_j\right]-\frac{\hbar^2 k^2 \Delta t}{24\pi^2}\left[\vbar^{l-1}_j+4\vbar^l_j+\vbar^{l+1}_j\right]\\*
&\qquad-\frac{\Delta t}{24\pi^2}\sideset{}{'}{\sum}_{n=j-N/2}^{N/2}\hat V_{n-j}\left[\vbar^{l-1}_n+4\vbar^l_n+\vbar^{l+1}_n\right]-\frac{\lambda_l}{2\pi}\vbar^l_j\, , & \text{for } j=0,\ldots,N/2-1,\\
0 &=\frac{i\hbar}{4\pi}\left[\hat v^{l+1}_j-\hat v^{l-1}_j\right]-\frac{\hbar^2 k^2 \Delta t}{24\pi^2}\left[\hat v^{l-1}_j+4\hat v^l_j+\hat v^{l+1}_j\right]\\*
&\qquad-\frac{\Delta t}{24\pi^2}\sideset{}{'}{\sum}_{n=j-N/2}^{N/2}\hat V_{j-n}\left[\hat v^{l-1}_n+4\hat v^l_n+\hat v^{l+1}_n\right]-\frac{\lambda_l}{2\pi}\hat v^l_j\, , & \text{for } j=0,\ldots,N/2-1,\\
0 &=\frac{i\hbar}{16\pi}\left[\vbar^{l-1}_{N/2}-\vbar^{l+1}_{N/2}\right]-\frac{\hbar^2 k^2 \Delta t}{96\pi^2}\left[\vbar^{l-1}_{N/2}+4\vbar^l_{N/2}+\vbar^{l+1}_{N/2}\right]\\*
&\qquad-\frac{\Delta t}{48\pi^2}\sideset{}{'}{\sum}_{n=0}^{N/2}\hat V_{n-N/2}\left[\vbar^{l-1}_n+4\vbar^l_n+\vbar^{l+1}_n\right]-\frac{\lambda_l}{2\pi}\vbar^l_{N/2}\, ,\\
0 &= \frac{i\hbar}{16\pi}\left[\hat v^{l+1}_{N/2} - \hat v^{l-1}_{N/2}\right]-\frac{\hbar^2 k^2 \Delta t}{96\pi^2}\left[\hat v^{l-1}_{N/2}+4\hat v^l_{N/2}+\hat v^{l+1}_{N/2}\right]\\*
&\qquad-\frac{\Delta t}{48\pi^2}\sideset{}{'}{\sum}_{n=0}^{N/2}\hat V_{N/2-n}\left[\hat v^{l-1}_n+4\hat v^l_n+\hat v^{l+1}_n\right]-\frac{\lambda_l}{2\pi}\hat v^l_{N/2}\, ,\\
1 &= \frac{1}{2\pi}\,\,\sideset{}{''}{\sum}_{k=-N/2}^{N/2}\hat v_k^l\vbar_k^l\, ,\\
0 &= \hat v^l_{-N/2} - \hat v^l_{N/2}\, ,\\
0 &= \vbar^l_{-N/2} - \vbar^l_{N/2}\, .
\end{align*}
This system of $(2N+3)$-equations, allow us to solve for $\{\hat v^{l+1}_k,\vbar^{l+1}_k\}_{k=-N/2}^{N/2}$ and $\lambda_l$ from initial data, $\{\hat v^{l-1}_k,\vbar^{l-1}_k\}_{k=-N/2}^{N/2}$ and $\{\hat v^l_k,\vbar^l_k\}_{k=-N/2}^{N/2}$. As such, this system of equations are an example of a spectral in space, second-order in time, {\bfi pseudospectral variational integrator} for the time-dependent Schr\"odinger equation. The expressions for the variational integrator for the time-independent Schr\"odinger equation, which has spectral accuracy in space, are given by
\begin{align*}
\hbar^2 k^2 \vbar_j
&=-\sideset{}{'}{\sum}_{n=-N/2}^{N/2+j}\hat V_{n-j}\vbar_n-\lambda\vbar_j\, ,
& \text{for } j=-N/2+1,\ldots,-1,\\
\hbar^2 k^2 \hat v_j
&=-\sideset{}{'}{\sum}_{n=-N/2}^{N/2+j}\hat V_{j-n}\hat v_n-\lambda\hat v_j\, ,
& \text{for } j=-N/2+1,\ldots,-1,\\
\hbar^2 k^2 \vbar_j
&=-\sideset{}{'}{\sum}_{n=j-N/2}^{N/2}\hat V_{n-j}\vbar_n-\lambda\vbar_j\, ,
& \text{for } j=0,\ldots,N/2-1,\\
\hbar^2 k^2 \hat v_j
&=-\sideset{}{'}{\sum}_{n=j-N/2}^{N/2}\hat V_{j-n}\hat v_n-\lambda\hat v_j\, ,
& \text{for } j=0,\ldots,N/2-1,\\
\frac{\hbar^2 k^2}{2} \vbar_{N/2}
&=-\sideset{}{'}{\sum}_{n=0}^{N/2}\hat V_{n-N/2}\vbar_n-\lambda\vbar_{N/2}\, ,\\
\frac{\hbar^2 k^2}{2} \hat v_{N/2}
&=-\sideset{}{'}{\sum}_{n=0}^{N/2}\hat V_{N/2-n}\hat v_n-\lambda\hat v_{N/2},\\
1&=\frac{1}{2\pi}\,\,\sideset{}{''}{\sum}_{k=-N/2}^{N/2}\hat v_k\vbar_k\, ,\\
\hat v^l_{-N/2} &= \hat v^l_{N/2}\, ,\\
\vbar^l_{-N/2} &= \vbar^l_{N/2}\, .
\end{align*}
As mentioned previously, it is possible generalize this approach to construct a fully spectral variational integrator in space-time, using Chebyshev polynomials to interpolate in time the coefficients of the discrete Fourier expansion used in the spatial interpolation. The computational cost of implementing such a scheme would be significantly higher, since this would require all the spatio-temporal degrees of freedom to be solved for simultaneously.
%%%%%%%%%%%%%%%%%%%%%%%%%%%%%%%%%%%%%%%%%%%%%%
\section{Conclusions and Future Work}
We have introduced the notion of a generalized Galerkin variational integrator, which is based on the idea of appropriately choosing a finite-dimensional approximation of the section of the configuration bundle, and approximating the action integral by a numerical quadrature scheme.

In contrast to standard variational methods, that are typically formulated in terms of interpolatory schemes parameterized by values of field variables at nodal and internal points, generalized Galerkin methods utilize function spaces that can be generated by arbitrary degrees of freedom. This allows the introduction of Lie group methods, and their symmetry reduction using discrete Euler--Poincar\'e reduction, as well as multiscale, and pseudospectral methods. Nonlinear approximation spaces allow the construction of spatio-temporally adaptive methods, which are better able to resolve shocks and other kinds of localized discontinuities in the solution.

It would be interesting to compare the performance of pseudospectral variational integrators with traditional pseudospectral schemes to see if any additional benefits arise from constructing pseudospectral schemes using a variational approach. More interesting still would be the comparison for fully spectral methods, since both variational and non-variational methods would achieve spectral accuracy, and it would make a particularly compelling case for variational integrators if their advantages  persist even when compared to numerical methods with spectral accuracy.

Most mesh adaptive methods use the principle of equipartitioning the error of the numerical scheme over the mesh elements to obtain moving mesh equations. These methods rely on \textit{a posteriori} error estimators that are related to the norm in which the accuracy of the numerical method is measured. While adaptive variational integrators exhibit an equipartitioning principle, in the sense that the discrete conjugate momentum associated with the horizontal variations are preserved from element to element in each connected component of the domain, it would be interesting to carefully explore the question of whether this can be understood as arising from error equipartitioning with respect to a geometrically motivated error estimator.

While we have only discussed the application of multiscale variational integrators to the case of ordinary differential equations, it would be natural to consider their generalizations to partial differential equations, whereby the multiscale shape functions are obtained through well-resolved solutions of the cell problem, as in the case with multiscale finite elements (see, for example, \cite{HoWu1999}). In general, short-term simulations at the fine scale can be used to construct appropriate shape functions to obtain generalized Galerkin variational integrators at a coarser level, through the use of principal orthogonal decomposition and balanced truncation, for example. This is consistent with the coarse-fine computational approach proposed in \cite{ThQiKe2000}, or the framework of heterogeneous multiscale methods as proposed in \cite{WeEn2003}.

A natural generalization would be to consider wavelet based variational integrators, as well as schemes based on conforming, hierarchical, adaptive refinement methods (CHARMS) introduced in \cite{GrKrSc2002} and further developed in \cite{KrTrZh2003}.

\bibliographystyle{plainnat}
\bibliography{umich_gvi}

\begin{thebibliography}{27}
\providecommand{\natexlab}[1]{#1}
\providecommand{\url}[1]{\texttt{#1}}
\expandafter\ifx\csname urlstyle\endcsname\relax
  \providecommand{\doi}[1]{doi: #1}\else
  \providecommand{\doi}{doi: \begingroup \urlstyle{rm}\Url}\fi

\bibitem[Baines(1995)]{Ba1995}
{M. J.} Baines.
\newblock \emph{Moving Finite Elements}.
\newblock Numerical Mathematics and Scientific Computation. Oxford University
  Press, 1995.

\bibitem[Chen and Hou(2003)]{ChHo2003}
Z.~Chen and {T. Y.} Hou.
\newblock A mixed multiscale finite element method for elliptic problems with
  oscillating coefficients.
\newblock \emph{Math. Comp.}, 72\penalty0 (242):\penalty0 541--576
  (electronic), 2003.

\bibitem[DeVore(1998)]{DeVore1998}
{R. A.} DeVore.
\newblock Nonlinear approximation.
\newblock In \emph{Acta Numerica}, volume~7, pages 51--150. Cambridge
  University Press, 1998.

\bibitem[E and Engquist(2003)]{WeEn2003}
W.~E and B.~Engquist.
\newblock The heterogeneous multiscale methods.
\newblock \emph{Commun. Math. Sci.}, 1\penalty0 (1):\penalty0 87--132, 2003.

\bibitem[Efendiev et~al.(2000)Efendiev, Hou, and Wu]{EfHoWu2000}
{Y. R.} Efendiev, {T. Y.} Hou, and {X. -H.} Wu.
\newblock Convergence of a nonconforming multiscale finite element method.
\newblock \emph{SIAM J. Numer. Anal.}, 37\penalty0 (3):\penalty0 888--910
  (electronic), 2000.

\bibitem[Grinspun et~al.(2002)Grinspun, Krysl, and Schr\"oder]{GrKrSc2002}
E.~Grinspun, P.~Krysl, and P.~Schr\"oder.
\newblock {CHARMS}: A simple framework for adaptive simulation.
\newblock \emph{ACM Transactions on Graphics (SIGGRAPH)}, 21\penalty0
  (21):\penalty0 281--290, July 2002.

\bibitem[Hou and Wu(1999)]{HoWu1999}
{T. Y.} Hou and {X. -H.} Wu.
\newblock A multiscale finite element method for {PDE}s with oscillatory
  coefficients.
\newblock In \emph{Numerical treatment of multi-scale problems (Kiel, 1997)},
  volume~70 of \emph{Notes Numer. Fluid Mech.}, pages 58--69. Vieweg, 1999.

\bibitem[Iserles(2003{\natexlab{a}})]{Is2003a}
A.~Iserles.
\newblock On the numerical quadrature of highly-oscillating integrals {I}:
  {F}ourier transforms.
\newblock Technical Report 2003/NA05, DAMTP, Cambridge, 2003{\natexlab{a}}.
\newblock (to appear in IMA J. Num. Anal.).

\bibitem[Iserles(2003{\natexlab{b}})]{Is2003b}
A.~Iserles.
\newblock On the numerical quadrature of highly-oscillating integrals {II}:
  {I}rregular oscillators.
\newblock Technical Report 2003/NA09, DAMTP, Cambridge, 2003{\natexlab{b}}.

\bibitem[Iserles(2004)]{Is2004}
A.~Iserles.
\newblock On the method of {N}eumann series for highly oscillatory equations.
\newblock Technical Report 2004/NA02, DAMTP, Cambridge, 2004.

\bibitem[Iserles and {N\o rsett}(2004)]{IsNo2004}
A.~Iserles and {S. P.} {N\o rsett}.
\newblock Efficient quadrature of highly oscillatory integrals using
  derivatives.
\newblock Technical Report 2004/NA03, DAMTP, Cambridge, 2004.

\bibitem[Iserles et~al.(2000)Iserles, Munthe-Kaas, {N\o rsett}, and
  Zanna]{IsMuNoZa2000}
A.~Iserles, H.~Munthe-Kaas, {S. P.} {N\o rsett}, and A.~Zanna.
\newblock Lie-group methods.
\newblock In \emph{Acta Numerica}, volume~9, pages 215--365. Cambridge
  University Press, 2000.

\bibitem[Jos{\'e} and Saletan(1998)]{JoSa1998}
{J. V.} Jos{\'e} and {E. J.} Saletan.
\newblock \emph{Classical Dynamics: {A} Contemporary Approach}.
\newblock Cambridge University Press, 1998.

\bibitem[Kane et~al.(1999)Kane, Marsden, and Ortiz]{KaMaOr1999}
C.~Kane, {J. E.} Marsden, and M.~Ortiz.
\newblock Symplectic-energy-momentum preserving variational integrators.
\newblock \emph{J. Math. Phys.}, 40\penalty0 (7):\penalty0 3353--3371, 1999.

\bibitem[Krysl et~al.(2003)Krysl, Trivedi, and Zhu]{KrTrZh2003}
P.~Krysl, A.~Trivedi, and B~. Zhu.
\newblock Object-oriented hierarchical mesh refinement with {CHARMS}.
\newblock \emph{Int. J. Numer. Meth. Eng.}, 2003.
\newblock (to appear).

\bibitem[Lall and West(2003)]{LaWe2003}
S.~Lall and M.~West.
\newblock Discrete variational mechanics and duality.
\newblock (in preparation), 2003.

\bibitem[Lew et~al.(2003)Lew, Marsden, Ortiz, and West]{LeMaOrWe2003}
A.~Lew, {J. E.} Marsden, M.~Ortiz, and M.~West.
\newblock Asynchronous variational integrators.
\newblock \emph{Arch. Ration. Mech. An.}, 167\penalty0 (2):\penalty0 85--146,
  2003.

\bibitem[Lew et~al.(2004)Lew, Marsden, Ortiz, and West]{LeMaOrWe2004}
A.~Lew, {J. E.} Marsden, M.~Ortiz, and M.~West.
\newblock Variational time integrators.
\newblock \emph{Int. J. Numer. Meth. Eng.}, 2004.
\newblock (to appear).

\bibitem[Marsden and West(2001)]{MaWe2001}
{J. E.} Marsden and M.~West.
\newblock Discrete mechanics and variational integrators.
\newblock In \emph{Acta Numerica}, volume~10, pages 317--514. Cambridge
  University Press, 2001.

\bibitem[Marsden et~al.(1998)Marsden, Patrick, and Shkoller]{MaPaSh1998}
{J. E.} Marsden, {G. W.} Patrick, and S.~Shkoller.
\newblock Multisymplectic geometry, variational integrators, and nonlinear
  {PDE}s.
\newblock \emph{Commun. Math. Phys.}, 199\penalty0 (2):\penalty0 351--395,
  1998.

\bibitem[Marsden et~al.(1999)Marsden, Pekarsky, and Shkoller]{MaPeSh1999}
{J. E.} Marsden, S.~Pekarsky, and S.~Shkoller.
\newblock Discrete {E}uler--{P}oincar\'e and {L}ie--{P}oisson equations.
\newblock \emph{Nonlinearity}, 12\penalty0 (6):\penalty0 1647--1662, 1999.

\bibitem[Marsden et~al.(2001)Marsden, Pekarsky, Shkoller, and
  West]{MaPeShWe2001}
{J. E.} Marsden, S.~Pekarsky, S.~Shkoller, and M.~West.
\newblock Variational methods, multisymplectic geometry and continuum
  mechanics.
\newblock \emph{J. Geom. Phys.}, 38\penalty0 (3-4):\penalty0 253--284, 2001.

\bibitem[Nocedal and Wright(1999)]{NoWr1999}
J.~Nocedal and {S. J.} Wright.
\newblock \emph{Numerical Optimization}.
\newblock Springer Series in Operations Research. Springer-Verlag, 1999.

\bibitem[Oliver et~al.(2004)Oliver, West, and Wulff]{OlWeWu2004}
M.~Oliver, M.~West, and C.~Wulff.
\newblock Approximate momentum conservation for spatial semidiscretizations of
  nonlinear wave equations.
\newblock \emph{Numer. Math.}, 2004.
\newblock (to appear).

\bibitem[Theodoropoulos et~al.(2000)Theodoropoulos, Qian, and
  Kevrekidis]{ThQiKe2000}
K.~Theodoropoulos, {Y. -H.} Qian, and {I. G.} Kevrekidis.
\newblock ``{C}oarse" stability and bifurcation analysis using timesteppers: a
  reaction diffusion example.
\newblock \emph{Proc. Natl. Acad. Sci.}, 97\penalty0 (18):\penalty0 9840--9843,
  2000.

\bibitem[Thoutireddy and Ortiz(2003)]{ThOr2003}
P.~Thoutireddy and M.~Ortiz.
\newblock A variational {$r$}-adaptation and shape-optimization method for
  finite-deformation elasticity.
\newblock \emph{Int. J. Numer. Meth. Eng.}, 2003.
\newblock (to appear).

\bibitem[Trefethen(2000)]{Trefethen2000}
{L. N.} Trefethen.
\newblock \emph{Spectral methods in {MATLAB}}.
\newblock Software, Environments, and Tools. Society for Industrial and Applied
  Mathematics (SIAM), 2000.

\end{thebibliography}

\end{document}